\documentclass[12pt,a4paper,twoside,english]{scrartcl}

\usepackage{graphicx}

\usepackage[a4paper,left=30mm,right=30mm,top=32mm,bottom=32mm]{geometry}
\usepackage{fleqn}
\usepackage{grffile}

\usepackage{amsmath}
\usepackage{amssymb}
\usepackage{amsfonts}
\usepackage{amsthm}
\usepackage{xcolor}

\usepackage[english]{babel}
\usepackage{inputenc}
\usepackage{nicefrac}
\usepackage{mathrsfs}
\usepackage{verbatim}

\newtheorem{thm}{Theorem}
\newtheorem{lem}[thm]{Lemma}
\newtheorem{prop}[thm]{Proposition}
\newtheorem{cor}[thm]{Corollary}

\newtheorem*{pf}{Proof}

\newtheorem*{pot2}{Proof of Theorem \ref{thm:ex2}}
\newtheorem*{pot3}{Proof of Theorem \ref{thm:ex}}
\newtheorem{ass}{Assumption}
\newtheorem{sys}{Problem}

\newcommand{\R}{\mathbb{R}}
\newcommand{\N}{\mathbb{N}}
\newcommand{\eps}{\varepsilon}
\newcommand{\calM}{\mathcal{M}}
\newcommand{\calL}{\mathcal{L}}
\newcommand{\embeds}{\hookrightarrow}
\DeclareMathOperator{\loc}{loc}
\DeclareMathOperator{\Id}{Id}

\begin{document}

\title{Analysis and numerical treatment of bulk-surface reaction-diffusion models of Gierer-Meinhardt type}

\author{Jan-Phillip~B\"acker\thanks{jan-phillip.baecker@tu-dortmund.de} $^,$\thanks{Institute of Applied Mathematics (LS III), TU Dortmund University, Vogelpothsweg 87, D-44227 Dortmund, Germany}\and
        Matthias~R\"oger\thanks{AG Biomathematik, TU Dortmund University, Vogelpothsweg 87, D-44227 Dortmund, Germany}\and
				Dmitri~Kuzmin$^\dagger$
}

\date{October 5, 2020}

\maketitle

\begin{abstract}
  We consider a Gierer-Meinhardt system on a surface coupled with a parabolic PDE in the bulk, the domain confined by this surface. Such a model was recently proposed and analyzed for two-dimensional bulk domains by Gomez, Ward and Wei ({\em SIAM J.~Appl.~Dyn.~Syst.~18}, 2019). We prove the well-posedness of the bulk-surface system in arbitrary space dimensions and show that solutions remain uniformly bounded in parabolic Hölder spaces for all times. The proof uses Schauders fixed point theorem and a splitting in a surface and a bulk part. We also solve a reduced system, corresponding to the assumption of a well mixed bulk solution, numerically. We use operator-splitting methods which combine a finite element discretization of the Laplace-Beltrami operator with a positivity-preserving treatment of the source and sink terms. The proposed methodology is based on the flux-corrected transport (FCT) paradigm. It constrains the space and time discretization of the reduced problem in a manner which provides positivity preservation, conservation of mass, and second-order accuracy in smooth regions. The results of numerical studies for the system  on a two-dimensional sphere demonstrate the occurrence of localized steady-state multispike pattern that have also been observed in one-dimensional models.
	
\end{abstract}

\section{Introduction}
The formation of patterns is a key property of many biological systems. Hence, insights into the nature of underlying organizing principles are of fundamental importance for various applications.
Several mechanisms have been identified that allow for the generation of spatial inhomogeneities in initially homogenous systems.
Turing observed in his seminal paper \cite{Turing1952} the possibility of diffusion induced instabilities that lead to symmetry breaking bifurcations and the formation of patterns.
Although diffusion is in general smoothing out spatial heterogeneities,strong variations of diffusivity in multi-component systems may lead to the enhancement of spatial structure. This typically requires a local self-activation of one species and a long-range antagonist.
Throughout the past decades numerous studies of pattern formation properties of reaction-diffusion systems have been presented.
The most relevant for this paper is a particular type of activator-inhibitor model proposed by Gierer and Meinhardt \cite{GM}. Their work has attracted considerable attention in the mathematics community, see for example \cite{NiTa86,MasTak1987,Wei99} and the review \cite{Ni98} for further references.
Most often the Gierer-Meinhardt system and its extensions are studied in open subsets of the Euclidean space, but corresponding systems are also well-motivated on lower-dimensional surfaces \cite{BeHa82}.

Recently, {\em bulk-surface} coupled Gierer-Meinhardt type systems have been analyzed in \cite{GoWaWe2018} for
an activator-inhibitor system on a one-dimensional membrane coupled to a diffusion process in the interior. The analysis presented in \cite{GoWaWe2018} is focused on membrane-bound spike patterns of this system.
Bulk-surface coupled reaction-diffusion systems appear in a wide range of applications in cell biology and other fields \cite{LeRa05,NGCR07,TLLW09,Miel13,RaeRoe14,GKRR16,ElRV17,BKMS17}. An increasing number of contributions have been made to the mathematical analysis of such models \cite{MoSh16,FeLT17,HaRo18,Diss20}.

The goal of this paper is to investigate the coupled bulk-surface reaction-diffusion system proposed in \cite{GoWaWe2018} in arbitrary space dimension and to provide a well-posedness analysis of the system and reliable numerical simulations.

In what follows, let $\Omega$ denote an open set in $\R^n$ with boundary $\Gamma=\partial\Omega$, representing the cell body and the cell membrane, respectively. Let $\nu$ denote the outer unit normal field of $\Omega$ and fix a time interval $(0,T)$. In the bulk-surface Gierer-Meinhardt system proposed in \cite{GoWaWe2018} the evolution of the concentration of two membrane bound proteins $u,v:\Gamma\times (0,T)\to\R$ and of a cytosolic protein $w:\Omega\times (0,T)\to\R$ are described. They satisfy the following reaction-diffusion system.

\begin{sys}[Bulk-surface Gierer-Meinhardt system]\label{sys:Ana}
Find a solution $(u,v,w)$ of
\begin{equation*}
\begin{split}
\partial_tu&=\varepsilon^2\Delta u-u+\frac{u^p}{v^q}+\sigma\hspace{100pt} \text{on}\ \Gamma\times(0,T),\\
\tau_s\partial_tv&=D_s\Delta v-(1+K)v+Kw+\varepsilon^{-1}\frac{u^r}{v^s}\hspace{37pt} \text{on}\ \Gamma\times(0,T),\\
\tau_b\partial_t w&=D_b\Delta w-w\hspace{137pt} \text{in}\ \Omega\times(0,T),\\
D_b\nabla w\cdot\nu&=Kv-Kw\hspace{141pt} \text{on} \ \Gamma\times(0,T),
\end{split}
\end{equation*}
subject to prescribed initial conditions.
\end{sys}
Precise conditions on the spatial domain, model parameters and initial conditions are fixed below. The equations for $u,v$ are an extension of the usual Gierer-Meinhardt system (the original system is obtained for $p=r=2$, $q=1$, $s=0$, $\sigma=K=0$). The variable $v$ is slowly diffusing compared to $u$, indicated by the small parameter $\eps>0$. The new contribution is the coupling to the additional variable $w$ that satisfies a diffusion equation in the bulk with a Robin-type boundary condition. In the second equation on the membrane, $w$ appears as a source term.

This system was analyzed and treated numerically by Gomez et al.~\cite{GoWaWe2018} in two dimensions. The authors presented a linear stability analysis of strongly localized structures for the spherical case and for the reduced model obtained in the limit $D_b\to \infty$. Additionally, they analyzed the occurrence of bifurcations and provided corresponding phase diagrams in parameter space. In this paper we show the existence of classical solutions to the full system Problem \ref{sys:Ana}. We use a fixed point argument and an operator splitting approach. Therefore the system is decomposed in a sub-system on the surface in the $u,v$ variables (for $w$ fixed), and the Robin boundary value problem in the bulk variable $w$ (for $v$ given). For both separate sub-systems we need suitable a-priori estimates and appropriate existence results. For the surface system we can suitably extend and adapt results and techniques from \cite{MasTak1987}, whereas parabolic theory of linear equations (see for example \cite{LaSoUr1968}) is available for the treatment of the bulk equation. The main focus is therefore on the derivation of suitable a-priori estimates for the two subproblems and the full system and an appropriate set-up for the Schauder fixed point theorem.

Our main results on the well-posedness of the bulk-surface system and properties of solutions are stated in Section \ref{sec:2}, proofs are given in Sections \ref{sub:Surf}-\ref{sub:Full}.

\bigskip
In addition to Problem \ref{sys:Ana} we investigate a reduced problem that arises in the limit of large bulk diffusivity $D_b\rightarrow \infty$ and leads to a well-mixed system where $w$ becomes spatially constant. This limit is motivated by the fact that cytosolic diffusion is typically much larger than lateral diffusion on the cell membrane.
\begin{sys}[Nonlocal surface Gierer-Meinhardt system]\label{sys:Num}
Find a solution $(u,v,w)$ of
\begin{equation*}
 \begin{split}
  \partial_tu&=\varepsilon^2\Delta_\Gamma u-u+\frac{u^p}{v^q}+\sigma \hspace{100pt} \text{on}\ \Gamma\times(0,T),\\
  \tau_s\partial_tv&=D_s\Delta_\Gamma v-(1+K)v+\frac{K}{|\Omega|}w+\varepsilon^{-1}\frac{u^r}{v^s}\hspace{32pt} \text{on}\ \Gamma\times(0,T),\\
  \tau_b\frac{\mathsf{d}}{\mathsf{d}t} w&=\left(1-K\frac{|\Gamma|}{|\Omega|}\right)w+K\int_\Gamma v \hspace{80pt}\text{in}\ (0,T),
 \end{split}
\end{equation*}
subject to prescribed initial conditions.
\end{sys}
This system has the advantage that all variables have their spatial domain of definition on the membrane only. It represents a kind of shadow system. Reductions to shadow systems have been analyzed intensively for two-variable reaction-diffusion systems \cite{Keen78,HaSa89,LiNi09,KaSu17,MHKS18}.  
The variable $w$ is determined by an ODE that is independent of $u$ but
depends on the total amount of $v$ on the membrane $\Gamma$.
Therefore it remains to solve a non-local membrane system in $u,v$. In Section \ref{sec:3} we discretize the reduced system in space using the continuous Galerkin finite element method. A tailor-made approximation of the source and sink terms is proposed to ensure positivity preservation for the solution of the fully discrete problem. Using first-order operator splitting \cite{HunVer2003}, row-sum mass lumping, and the classical Patankar method \cite{BuDeMe2003,Patank1980} for semi-implicit time integration of nonlinear source terms, we construct a first-order positivity-preserving scheme for the coupled problem. Its second-order counterpart is designed using the second-order accurate Strang splitting \cite{Strang1968}, an adaptive combination \cite{Gruell2018} of the Patankar method with a second-order strong stability preserving (SSP) Runge-Kutta time discretization of the source terms, and selective mass lumping based on the flux-corrected transport (FCT) methodology \cite{Kuzmin2011,Zalesak1979}. In Section \ref{sec:4} simulations are performed for different choices of the initial data and different values of $K$. The results of these simulations are compared to the results of Gomez et al. \cite{GoWaWe2018}.

In Section \ref{sec:4} we will present various simulations that show the occurence and preservation of strongly localized structures in this system. Instead of  a spatially homogenous solution, long-time simulations produce different solutions with spike patterns. We can observe the occurence of patterns with a different number of spikes. The most stable states are those with one or two spikes and configuration with symmetrically distributed spikes. This pattern formation behavior is consistent with the results of simulations on the circle and the findings of Gomez et al.~\cite{GoWaWe2018}.

\subsection{Notation}
Fix an open set $\Omega\subset\mathbb{R}^n$ with smooth boundary  $\Gamma=\partial\Omega$ and $0<T\leq\infty$. Denote the corresponding space-time cylinder by $\Omega_T$ and its lateral boundary by  $\Gamma_T=\Gamma\times (0,T)$.

For a subset $A$ in $\mathbb{R}^n$, an interval $I\subset [0,\infty)$ and $k\in\N_0$ we define
\begin{equation*}
 \begin{split}
  C^{2k,k}(A\times I) &:=\big\{u:A\times I\to\R\,:\,\partial_t^lD^{\gamma}u \in C^0(A\times I)\\
  &\qquad\qquad\qquad \text{for all}\ 2l+|\gamma|\leq 2k,\ l\in\mathbb{N}_0, \gamma\in(\mathbb{N}_0)^n\big\}.
 \end{split}
\end{equation*}

We next briefly introduce parabolic Hölder spaces. We define the parabolic distance in $\R^n\times\R$ by
\begin{equation*}
 d\big((x,t),(y,s)\big) = \big(|x-y|^2+|t-s|\big)^{\frac{1}{2}}.
\end{equation*}
For $U=\Omega$ or $U=\Gamma$, an open interval $I\subset (0,\infty)$, and $0<\lambda\leq 1$ we define for $u:U_T\to\R$ the parabolic Hölder constants
\begin{equation*}
 [u]_{U_T,\lambda} := \sup_{(x,t)\neq (y,s)\in U_T} \frac{|u(x,t)-u(y,s)|}{d\big((x,t),(y,s)\big)}
\end{equation*}
and for $k\in \N_0$ the parabolic Hölder spaces by
\begin{equation*}
 C^{2k+\lambda,k+\frac{\lambda}{2}}(U_T):=\{u:U_T\to\R\,:\, \|u\|_{C^{2k+\lambda,k+\frac{\lambda}{2}}(U_T)}<\infty\},
\end{equation*}
where
\begin{equation*}
 \|u\|_{C^{2k+\lambda,k+\frac{\lambda}{2}}(U_T)} :=
 \sum_{2l+|\gamma|\leq 2k} \|\partial_t^lD^{\gamma} u\|_{L^\infty(U_T)} + [\partial_t^lD^{\gamma}u]_{U_T,\lambda}.
\end{equation*}
We remark that for $u\in C^{2k+\lambda,k+\frac{\lambda}{2}}(U_T)$ all derivatives $\partial_t^lD^{\gamma} u$ with $2l+|\gamma|\leq 2k$ can be continuously extended to $\overline{U_T}$ and $C^{2k+\lambda,k+\frac{\lambda}{2}}(U_T)\embeds C^{2k,k}(\overline{U_T})$ is compact.

For $U_T\subset A\subset \overline{U_T}$ we set
\begin{align*}
 C^{2k+\lambda,k+\frac{\lambda}{2}}_{\loc}(A) := \Big\{ u\in C^{2k,k}(U_T) :
 \text{ for all }(x,t)\in A \text{ exists }\rho>0\text{ such that }\\
  \quad
  u\in C^{2k+\lambda,k+\frac{\lambda}{2}}\big(U_T\cap \big\{d\big(\cdot,(x,t)\big)<\rho\big\}\big) \Big\}.
\end{align*}

We use the shorthand notation $W^{2,1}_r(\Gamma_T),\ 1\leq r\leq\infty$ for the parabolic Sobolev space $W^{1,r}(0,T;L^r(\Gamma))\cap L^r(0,T;W^{2,r}(\Gamma))$.

We call a solution $(u,v,w)$ of Problem \ref{sys:Ana} a classical solution, if $u,v$ belong to the class $C^{2,1}(\Gamma_T)\cap C^0([0,T)\times\Gamma)$ and $w$ belongs to the class $C^{2,1}(\Omega_T)\cap C^0((0,T);C^1(\overline{\Omega}))\cap C^0([0,T)\times\Omega)$. We call a solution $(u,v,w)$ a positive solution, if $u$ and $w$ are non-negative and $v$ is positive. In an analogous way we define classical solutions and positive solutions for the membrane sub-system Problem \ref{sys:Surf} and the bulk equation Problem \eqref{sys:Bulk} introduced below.

We will use references to standard regularity results for parabolic problems on subsets of $\mathbb{R}^n$ also for the case in which the spatial domain is a submanifold. This can be justified by a partition of unity subordinate to parametrized subsets of the manifold, see for example \cite{HaRo18}.

In the statements and calculations below $C$ typically denotes a generic constant independent of the time $T$, while $\Lambda$ denotes a fixed constant independent of $T$.

\section{Well-posedness of the bulk-surface Gierer-Meinhardt model}
\label{sec:2}
In this section we formulate our main analytical results for the fully coupled system Problem \ref{sys:Ana}. The following assumptions are used.

\begin{ass}\label{ass:Omega}
Let $\Omega\subset\mathbb{R}^n$ be an open, bounded and connected subset with smooth boundary, $\Gamma=\partial\Omega$ and $\nu $ be the outer unit normal of $\Omega$. Assume, that $\tau_s,\tau_b,D_s,D_b,K>0$, $\varepsilon,\sigma>0$ hold and $p>1,q>0,r>0,s\geq0$ satisfy
\begin{equation*}
 0<\frac{p-1}{r}<\frac{q}{s+1},\qquad \frac{p-1}{r}<\frac{2}{n+1}.
\end{equation*}
Let $u_0,v_0\in C^{2+\lambda}(\Gamma)$, $w_0\in C^{2+\lambda}(\Omega)$, $\min_\Gamma u_0\geq 0$, $\min_\Gamma v_0,\min_{\overline{\Omega}}w_0>0$ hold for some $\lambda\in(0,1)$. We prescribe the compatibility condition
\begin{equation*}
 D_b\nabla w_0\cdot\nu = K(v_0-w_0)\quad\text{ on }\Gamma.
\end{equation*}

We choose $\Lambda_0$ such that
\begin{equation*}
 \max\Big\{\|u\|_{C^{2+\lambda}(\Gamma)}, \|u\|_{C^{2+\lambda}(\Gamma)},
 \|w\|_{C^{2+\lambda}(\Omega)}\Big\}\,\leq\, \Lambda_0.
\end{equation*}
\end{ass}

Our main theorem presents the global-in-time existence, uniqueness and regularity of classical solutions of Problem \ref{sys:Ana}.

\begin{thm}[Existence, regularity and uniqueness of solutions]\label{thm:ex}
Let Assumption \ref{ass:Omega} hold. Then there exists a unique, positive classical solution $(u,v,w)$ of Problem \ref{sys:Ana} with initial data
\begin{align}
 u(\cdot,0) &= u_0,\quad v(\cdot,0) = v_0 \qquad
\text{ on }\Gamma  \label{eq:ini-uv},\\
w(\cdot,0) &= w_0 \qquad
\text{ in }\Omega.  \label{eq:ini-w}
\end{align}
This solution exists globally in time and satisfies
$u,v\in C^{2+\lambda,\frac{2+\lambda}{2}}(\Gamma\times [0,\infty))$, $w\in C^{2+\lambda,\frac{2+\lambda}{2}}(\Omega\times [0,\infty))$.
\end{thm}

Note that this result in particular implies that $u(\cdot,t),v(\cdot,t)$ and $w(\cdot,t)$ are uniformly bounded in $C^{2+\gamma}(\Omega)$ and $C^{2+\gamma}(\Gamma)$, respectively, independently of $t>0$.

This theorem will be proved in Sections \ref{sub:Surf}-\ref{sub:Full} by a fixed point argument making use of operator splitting. Therefore we will decompose the problem in a surface system and a boundary value problem in the bulk. These subproblems will be analyzed separately before considering the full system.

\begin{sys}[Surface system]\label{sys:Surf}
Let Assumption \ref{ass:Omega} hold, and let $0<T<\infty$ and a nonnegative function $\tilde{w}$ be given with
\begin{equation}
 \tilde{w}\in C^{\lambda,\frac{\lambda}{2}}(\Omega\times (0,T)).
 \label{eq:wtilde}
\end{equation}
Find a solution $(u,v)$ of
\begin{equation*}
\begin{split}
\partial_tu&=\varepsilon^2\Delta_\Gamma u-u+\frac{u^p}{v^q}+\sigma \hspace{100pt} \text{on}\ \Gamma\times(0,T),\\
\tau_s\partial_tv&=D_s\Delta_\Gamma v-(1+K)v+K\tilde{w}+\varepsilon^{-1}\frac{u^r}{v^s} \hspace{37pt} \text{on}\ \Gamma\times(0,T),
\end{split}
\end{equation*}
that in addition satisfies the initial condition \eqref{eq:ini-uv}.
\end{sys}

In Section \ref{sub:Surf} we prove the existence of solutions to this system and derive appropriate estimates for the solution.

Afterward, in Section \ref{sub:Bulk}, we will prove the existence of solutions of the following Robin boundary value problem in the bulk.
\begin{sys}[Bulk equation]\label{sys:Bulk}
Let Assumption \ref{ass:Omega} hold, and let $0<T<\infty$ and a positive function $\tilde{v}$ be given with
\begin{equation}
 \tilde{v}\in C^{1+\lambda,\frac{1+\lambda}{2}}(\Gamma\times(0,T)),\qquad
 \tilde v(\cdot,0)=v_0\quad\text{ on }\Gamma.
 \label{eq:vtilde}
\end{equation}
Find a solution $w$ of
\begin{equation*}
\begin{split}
\tau_b\partial_t w&=D_b\Delta w-w \hspace{37pt}\text{in}\ \Omega\times(0,T),\\
D_b\nabla w\cdot\nu&=K\tilde{v}-Kw\hspace{41pt}\text{on} \ \Gamma\times(0,T),
\end{split}
\end{equation*}
that in addition satisfies the initial condition \eqref{eq:ini-w}.
\end{sys}

In Section \ref{sub:Full} we will finally prove Theorem \ref{thm:ex}. We will introduce a suitable iteration map and use a fixed point argument. The main task is to verify the criteria for the Schauder fixed point theorem and to prove  global-in-time estimates.

\section{Surface system}\label{sub:Surf}
In this section we deal with the $u,v$ system on the surface and prove the following existence theorem.
\begin{thm}\label{thm:ex2}
Let $0<T<\infty$ be arbitrary and let Assumption \ref{ass:Omega} hold. Then for any $\tilde w\geq 0$ with \eqref{eq:wtilde} there exists a unique, non-negative solution $(u,v)\in (C^{2+\lambda,\frac{2+\lambda}{2}}(\Gamma_T))^2$ of Problem \ref{sys:Surf} with initial condition \eqref{eq:ini-uv}.

The solution $(u,v)$ satisfies for all $T<\infty$
\begin{equation*}
 \|u\|_{C^{2+\lambda,\frac{2+\lambda}{2}}(\Gamma_T)} +
 \|v\|_{C^{2+\lambda,\frac{2+\lambda}{2}}(\Gamma_T)}
 \leq C(\lambda,n,\Omega,T) \big(1+ \Lambda_0 + \|\tilde w\|_{C^{\lambda,\frac{\lambda}{2}}(\Omega_T)}\big),
\end{equation*}
where $C(\lambda,n,\Omega,T)$ remains bounded with $T\downarrow 0$.

The solution operator $\calM:C^{\lambda,\frac{\lambda}{2}}(\Omega_T)\to C^{2+\lambda,\frac{2+\lambda}{2}}(\Gamma_T)^2$ that maps $\tilde w$ to $(u,v)$ is continuous.
\end{thm}

Corresponding results in the case of a spatial domain given by an open subset of $\R^n$ and for $K=0$ have been shown in \cite{MasTak1987}. For the bulk-surface system we have in particular to take care of the additional term $K\tilde{w}$ in the PDE for $v$.
We follow the arguments in \cite{MasTak1987} and adapt their approach to the present bulk-surface system.
The proof uses a series of a-priori estimates. We start with lower bounds for the solutions and then prove in a bootstrapping procedure and with a number of intermediate steps the boundedness in parabolic Hölder spaces.

Throughout the whole section we fix an arbitrary $0<T<\infty$ and require that Assumption \ref{ass:Omega} and \eqref{eq:wtilde} hold. For the following lemmas and corollaries that deal with a priori estimates we assume a classical solution $(u,v)$ of Problem \ref{sys:Surf} on $[0,T)$.

\begin{lem}[cf \cite{MasTak1987}]\label{lem:Unten2}
There exist constants $m_u\geq 0,\ m_v>0$ independent of $T$ and $\tilde{w}$ such that
\begin{equation*}
 \inf_{\Gamma_T}u\geq m_u,\qquad
 \inf_{\Gamma_T}v\geq m_v.
\end{equation*}
\end{lem}

The proof of this lemma can be reduced to the proof of Lemma 2.1 in \cite{MasTak1987}. Since the sign of $\tilde{w}$ is positive, we can estimate the solution from below by the solution corresponding to $\tilde w=0$, for which the results in \cite{MasTak1987} apply.

We next obtain an estimate for the $L^\infty(0,T;L^k(\Gamma))$-norm of $u$.

\begin{lem}[cf \cite{MasTak1987}]\label{lem:Lk1}
Let $k\geq 1$ be arbitrary. Then there exists a positive constant $C_k$ independent of $T$ and $\tilde{w}$ such that
\begin{equation}
 \|u(\cdot,t)\|_{L^k(\Gamma)}\leq C_k(1+\Lambda_0)
 \label{eq:uLk}
\end{equation}
holds for $0<t<T$.
\end{lem}
\begin{pf}
The proof follows Lemma 2.3 and Lemma 2.4 in \cite{MasTak1987}. The first step is a control of the functions $v_\alpha := u^rv^{-(\alpha+s+1)}$, $\alpha>0$. It turns out that for any $0<t<T$ the function $\hat v$, defined by $\hat v(x,\tau):= e^{-k(t-\tau)}v(x,\tau)^{-\alpha}$ satisfies
\begin{align*}
 \partial_\tau \hat v - \frac{D_s}{\tau_s}\Delta\hat v
 &\leq k\hat v - \frac{\alpha(1+K)}{\tau_s}\hat v - \frac{\alpha K}{\tau_s}\frac{\hat v}{v}\tilde w -
 \frac{\alpha}{\eps\tau_s} e^{-k(t-\tau)}v_\alpha\\
 &\leq C(k,K,\alpha,\tau_s)e^{-k(t-\tau)}m_v^{-\alpha} -
 \frac{\alpha}{\eps\tau_s} e^{-k(t-\tau)}v_\alpha.
\end{align*}
Note that for $\tilde w$ we only have used the nonnegativity of $\tilde w$.
Integration over $\Gamma\times (0,t)$ yields the estimate
\begin{align}
 \frac{\alpha}{\eps\tau_s}\int_0^t e^{-k(t-\tau)}\int_\Gamma v_\alpha(x,\tau)\,dx\,d\tau
 &\leq C(k,K,\alpha,\tau_s,m_v,\Gamma) +e^{-kt}m_v^{-\alpha}|\Gamma| \notag\\
 &\leq C(k,K,\alpha,\tau_s,m_v,\Gamma).
 \label{eq:estalpha}
\end{align}
We next prove the required bound for $u$. This argument only uses the PDE for $u$ and is therefore the same as in \cite{MasTak1987}. Regarding $v$, only the estimate \eqref{eq:estalpha} is used. It is the proof of this lemma, where the conditions on the parameters $p,q,r,s$ from Assumption \ref{ass:Omega} are used. We refer to \cite[Lemma 2.4]{MasTak1987} for the details.
\hfill $\square$
\end{pf}

In the next step we show that a corresponding bound is  satisfied by $v$.
\begin{lem}\label{lem:Lk3}
For any $k\geq 1$ there exists a positive constant $C_k$ independent of $\tilde{w}$ and $T$, such that
\begin{equation*}
 \sup_{0<t<T}\|v(\cdot,t)\|_{L^k(\Gamma)}\leq C_k\big(1+\Lambda_0+\|\tilde w\|_{L^\infty(0,T;L^k(\Gamma))}\big).
\end{equation*}
\end{lem}

\begin{pf}
An integration by parts and the PDE for $v$ yield
\begin{equation*}
\begin{split}
\tau_s\frac{\mathsf{d}}{\mathsf{d} t}\|v\|_{L^k(\Gamma)}^k&=\int_\Gamma\tau_s\partial_t(v^k)=k\int_\Gamma v^{k-1}\tau_s\partial_t v\\&=k\int_\Gamma v^{k-1} \left(D_s\Delta_\Gamma v-(1+K)v+K\tilde{w}+\varepsilon^{-1}\frac{u^r}{v^s}\right)\\&=-k(1+K)\|v\|_{L^k(\Gamma)}^k-k(k-1)D_s\int_\Gamma v^{k-2}|\nabla_\Gamma v|^2\\&\quad+k\int_\Gamma K\tilde{w}v^{k-1}+k\varepsilon^{-1}\int_\Gamma\frac{u^r}{v^s}v^{k-1}.
\end{split}
\end{equation*}
Applying the Hölder and Young inequalities we derive
\begin{equation*}
 \begin{split}
  \tau_s\frac{\mathsf{d}}{\mathsf{d} t}\|v\|_{L^k(\Gamma)}^k\leq&-k(1+K)\|v\|_{L^k(\Gamma)}^k-k(k-1)D_s\int_\Gamma v^{k-2}|\nabla_\Gamma v|^2\\&+k\int_\Gamma K\delta_1 v^{k}+
  \int_\Gamma C_k(\delta_1)\tilde{w}^k\\
  &+k\varepsilon^{-1}\delta_2\int_\Gamma v^{k}
  +K\int_\Gamma\varepsilon^{-1}C_k(\delta_2)\frac{u^{kr}}{v^{ks}},
 \end{split}
\end{equation*}
with $C_k(\delta_1),C_k(\delta_2)>0$.
Since $v$ is bounded from below by $m_v$ we choose $\delta_1+ \varepsilon^{-1}\delta_2 K^{-1}=1$ and find that
\begin{equation*}
\tau_s\frac{\mathsf{d}}{\mathsf{d} t}\|v\|_{L^k(\Gamma)}^k\leq-k\|v\|_{L^k(\Gamma)}^k-k(k-1)D_s\int_\Gamma v^{k-2}|\nabla_\Gamma v|^2+\int_\Gamma C_1\tilde{w}^k+\int_\Gamma C_2u^{kr}
\end{equation*}
holds. It follows
\begin{equation*}
\tau_s\frac{\mathsf{d}}{\mathsf{d} t}\|v\|_{L^k(\Gamma)}^k\leq-k\|v\|_{L^k(\Gamma)}^k + C_1\|\tilde{w}\|_{L^k(\Gamma)}^{k}+C_2\|u\|_{L^{kr}(\Gamma)}^{kr}.
\end{equation*}
This gives
\begin{align*}
 &\|v(\cdot,t)\|_{L^k(\Gamma)}^k \\
 \leq\, &
 e^{\frac{-kt}{\tau_s}}\|v_0\|_{L^k(\Gamma)}^k +\frac{C_1+C_2}{\tau_s}\int_0^t e^{\frac{-k(t-\tau)}{\tau_s}}\Big(\|\tilde{w}(\cdot,\tau)\|_{L^k(\Gamma)}^{k}+\|u(\cdot,\tau)\|_{L^{kr}(\Gamma)}^{kr}\Big)\,d\tau\\
 \leq\,& e^{\frac{-kt}{\tau_s}}\|v_0\|_{L^k(\Gamma)}^k +C(\tau_s,k)\Big(\|\tilde{w}\|_{L^\infty(0,T;L^k(\Gamma))}^{k}+\|u\|_{L^\infty(0,T;L^{kr}(\Gamma))}^{kr}\Big).
\end{align*}
Using Lemma \ref{lem:Lk1} this implies the desired estimate.
\hfill$\square$
\end{pf}

Now we can also show higher regularity for the solutions.

\begin{lem}\label{lem:W1N}
Let $N\geq n-1$ be arbitrary. Then there exists a positive constant $C$ independent of $T$ and $\tilde{w}$, such that
\begin{align}
 \|\nabla u(\cdot,t)\|_{L^N(\Gamma)}
 &\leq C\big(1+\Lambda_0\big),
 \label{eq:uW1N}\\
 \|\nabla v(\cdot,t)\|_{L^N(\Gamma)}
 &\leq C\big(1+\Lambda_0+ \|\tilde{w}\|_{L^\infty(0,T;L^N(\Gamma))}\big)
 \label{eq:vW1N}
\end{align}
hold for all $0<t<T$.
\end{lem}

\begin{pf}
 The proof of this lemma is analogous to the proof of Lemma 2.7 in \cite{MasTak1987}. One uses estimates for the linear semigroup associated to the operator $A:=-\Delta + \Id$ on $\Gamma$, in particular the estimate
 \begin{equation*}
  \|e^{-tA}\|_{\calL(L^N(\Gamma);W^{1,N}(\Gamma))} \leq C(1+t^{-\frac{1}{2}})e^{-t}.
 \end{equation*}

 The arguments to derive the bounds for $u$ and $v$ are similar. We therefore consider the $v$ equation, since here a dependence on $\tilde w$ appears. The previous estimates then imply that
 \begin{align*}
  \|\nabla v\|_{L^N(\Gamma)} \leq Ce^{-t}\|v_0\|_{W^{2,N}(\Gamma)} + C(\tau_s,D_s)\int_0^t e^{-(t-\tau)}(1+|t-\tau|^{1/2})\|f_2(\cdot,\tau)\|_{L^N(\Gamma)}\,d\tau,
 \end{align*}
 where $f_2=-K\tilde w+\eps^{-1}\frac{u^r}{v^s}$. We then can use the
 lower bound $v\geq m_v$ and the $L^\infty(0,T;L^N(\Gamma))$ estimates for $\tilde w$ and $u$ to control the right-hand side and deduce \eqref{eq:vW1N}.
\hfill$\square$
\end{pf}

\begin{cor}\label{lem:Linf}
There exists a positive constant $C$ independent of $T$ and $\tilde{w}$ such that
\begin{align*}
 \|u(\cdot,t)\|_{L^{\infty}(\Gamma)}+
 \|v(\cdot,t)\|_{L^{\infty}(\Gamma)} \leq C\big(1+\|\tilde{w}\|_{L^\infty(\Omega_T)}\big)
\end{align*}
holds for all $0<t<T$.
\end{cor}

This corollary follows directly from Lemma \ref{lem:W1N} and Sobolev embeddings.

\begin{lem}\label{lem:Hoel1}
Let $0<T<\infty$. Then we have $u,v\in C^{2+\lambda,\frac{2+\lambda}{2}}(\Gamma_T)$ with the estimates
\begin{align*}
\|u\|_{C^{2+\lambda,\frac{2+\lambda}{2}}(\Gamma_T)}+
\|v\|_{C^{2+\lambda,\frac{2+\lambda}{2}}(\Gamma_T)}\leq
C(\lambda,n,\Omega,T) (1+ \Lambda_0 + \|\tilde w\|_{C^{\lambda,\frac{\lambda}{2}}(\Omega_T)}),
\end{align*}
where $C(\lambda,n,\Omega,T)$ remains bounded with $T\downarrow 0$.
\end{lem}
\begin{pf}
We write the PDE for $u$ in the form
\begin{equation*}
 \partial_t u-\varepsilon^2\Delta_\Gamma u=f_1
\end{equation*}
with $f_1=\sigma+\frac{u^p}{v^q}$. Since $u$ is in $L^k(0,T;L^k(\Gamma))$ for all $k$ and $v$ is greater than a positive number, $f_1$ belongs to $L^N(0,T;L^N(\Gamma))$ for any $1\leq N<\infty$. We fix $N>n$.
Regularity theory \cite[Theorem IV.9.1]{LaSoUr1968} implies that $u\in W^{2,1}_N(\Gamma_T)$ with
\begin{align*}
 \|u\|_{W^{2,1}_N(\Gamma_T)} &\leq
 C(\eps,N,n,\Omega,T)\Big(\|u_0\|_{W^{2,N}(\Gamma)}+\|f_1\|_{L^N(\Gamma_T)}\Big) \notag\\
 &\leq C(m_v,\eps,N,n,\Omega,T)(1+\Lambda_0),
\end{align*}
where in the last inequality we have used \eqref{eq:uLk} and where $C(m_v,\eps,N,n,\Omega,T)$ remains bounded with $T\downarrow 0$.

A similar estimate holds for $v$, where a dependence on $\|\tilde w\|_{L^N(\Gamma_T)}$ appears on the right-hand side. We therefore obtain $v\in W^{2,1}_N(\Gamma_T)$ with
\begin{equation*}
 \|v\|_{W^{2,1}_N(\Gamma_T)} \leq C(m_v,\eps,N,n,\Omega,T)\big(1+\Lambda_0+\|\tilde w\|_{C^{\lambda,\frac{\lambda}{2}}(\Omega_T)}\big).
\end{equation*}

Next consider arbitrary $0<s<t<T$.
We observe that
\begin{equation*}
 \|u(t)-u(s)\|_{L^N(\Gamma)} \leq \int_s^t \|\partial_t u(\tau)\|_{L^N(\Gamma)}\,d\tau 
 \leq \|\partial_t u(\tau)\|_{L^N(\Gamma_T)}(t-s)^{1-\frac{1}{N}},
\end{equation*}
and obtain from the Gagliardo-Nirenberg inequality for any $0<\theta<1$ that
\begin{align*}
 \|u(t)-u(s)\|_{W^{\theta,N}(\Gamma)} &\leq C\|u(t)-u(s)\|_{W^{1,N}(\Gamma)}^\theta\|u(t)-u(s)\|_{L^N(\Gamma)}^{1-\theta}
 \notag\\
 &\leq C(N,\theta)(1+\Lambda_0^\theta)\| u\|_{W^{1,N}(0,T;L^N(\Gamma))}^{1-\theta}|t-s|^{(1-\frac{1}{N})(1-\theta)}
 \notag\\
 &\leq C(m_v,\eps,N,n,\Omega,T,\theta)(1+\Lambda_0)|t-s|^{(1-\frac{1}{N})(1-\theta)},
\end{align*}
where we have used \eqref{eq:uW1N} and where $C(m_v,\eps,N,n,\Omega,T,\theta)$ remains bounded with $T\downarrow 0$.
By Sobolev-Morrey inequality this yields that for $\theta-\frac{2}{N}>0$
\begin{align*}
 \|u(t)-u(s)\|_{C^0(\Gamma)} &\leq C(m_v,\eps,N,n,\Omega,T,\theta)(1+\Lambda_0)|t-s|^{(1-\frac{1}{N})(1-\theta)}.
\end{align*}
Since also for any $0<\gamma<1-\frac{2}{N}$
\begin{equation*}
 \|u\|_{L^\infty(0,T;C^{0,\gamma}(\Gamma))} \leq
 C(\Omega,\gamma,N)\|u\|_{L^\infty(0,T;W^{1,N}(\Gamma))}
 \leq C(m_v,\eps,N,n,\Omega,T,\gamma)(1+\Lambda_0)
\end{equation*}
by \eqref{eq:vW1N}, and since $N<\infty$ is arbitrary we deduce that $u\in C^{\gamma,\gamma/2}(\Gamma_T)$ for all $0<\gamma<1$ with
\begin{align*}
 \|u\|_{C^{\gamma,\gamma/2}(\Gamma_T)} &\leq C(m_v,\eps,N,n,\Omega,T,\gamma)(1+\Lambda_0),
\end{align*}
where $C(m_v,\eps,N,n,\Omega,T,\gamma)$ remains bounded with $T\downarrow 0$. Similarly we deduce
\begin{align*}
 \|v\|_{C^{\gamma,\gamma/2}(\Gamma_T)} &\leq C(m_v,\eps,N,n,\Omega,T,\gamma)(1+\Lambda_0+\|\tilde w\|_{C^{\lambda,\frac{\lambda}{2}}(\Omega_T)}).
\end{align*}
Since $v\geq m_v>0$ and the nonlinearities on the right-hand side of the PDEs for $u,v$ are locally Lipschitz continuous functions of $u,v$, we deduce from Schauder estimates \cite[Theorem IV.5.2]{LaSoUr1968} that
\begin{align*}
 \|u,v\|_{C^{2+\lambda,1+\lambda/2}(\Gamma_T)} &\leq C(m_v,\eps,N,n,\Omega,T,\lambda)(1+\Lambda_0+\|\tilde w\|_{C^{\lambda,\frac{\lambda}{2}}(\Omega_T)}),
\end{align*}
where $C(m_v,\eps,N,n,\Omega,T,\lambda)$ remains bounded with $T\downarrow 0$ (see Theorem IV.5.4 and the proof of Theorem IV.5.2 in \cite{LaSoUr1968}).
\hfill$\square$
\end{pf}

The previous results eventually yield Theorem \ref {thm:ex2}.

\begin{pot2}
Lemma \ref{lem:Unten2} gives the lower bounds for the solutions, while Corollary \ref{lem:Linf} gives the upper bounds. The proof of existence and uniqueness follows as in Masuda and Takahashi \cite{MasTak1987} from local existence results for reaction-diffusion systems \cite{Roth84} and the preceding a-priori estimates.
The regularity statement follows from Lemma \ref{lem:Hoel1}.

To prove that $\calM$ is continuous we consider $\tilde w_1,\tilde w_2 \in C^{\lambda,\frac{\lambda}{2}}(\Omega_T)$ and let $(u_j,v_j)=\calM(\tilde w_j)$, $j=1,2$. Since the right-hand sides of the PDEs for $u_j,v_j$ are Lipschitz continuously depending on $u_j,v_j$ and linear in $\tilde w_j$ we deduce from Schauder estimates for systems \cite[Theorem VII.5.2]{LaSoUr1968}
\begin{align*}
 \|u_1-u_2\|_{C^{2+\lambda,\frac{2+\lambda}{2}}(\Gamma_T)}+
 \|v_1-v_2\|_{C^{2+\lambda,\frac{2+\lambda}{2}}(\Gamma_T)}
 &\leq C\|\tilde w_1-\tilde w_2\|_{C^{\lambda,\lambda/2}(\Gamma_T)}.
\end{align*}
\end{pot2}

\section{Bulk system}\label{sub:Bulk}
The existence theorem for the initial boundary problem in the bulk as formulated in Problem \ref{sys:Bulk} follows readily from the theory of linear parabolic equations.

\begin{thm}\label{thm:ex1}
Let $0<T<\infty$ be arbitrary and let Assumption \ref{ass:Omega} hold.
Then for any $\tilde v$ as in \eqref{eq:vtilde} there exists a unique positive solution $w\in C^{2+\lambda,\frac{2+\lambda}{2}}(\Omega_T)$ of Problem \ref{sys:Bulk}. Moreover, we have
\begin{equation}
 \|w\|_{C^{2+\lambda,\frac{2+\lambda}{2}}(\Omega_T)} \leq C(\lambda,n,\Omega,T)\big(1+\Lambda_0+\|\tilde v\|_{C^{1+\lambda, \frac{1+\lambda}{2}}(\Gamma_T)}\big),
 \label{eq:C2lambdaw}
\end{equation}
where $C(\lambda,n,\Omega,T)$ remains bounded with $T\downarrow 0$.

The solution operator $\calL:C^{1+\lambda,\frac{1+\lambda}{2}}(\Omega_T)\to C^{2+\lambda,\frac{2+\lambda}{2}}(\Omega_T)$ that maps $\tilde v$ to $w$ is continuous.
\end{thm}

\begin{pf}
The existence, uniqueness and the estimate \eqref{eq:C2lambdaw} follow from \cite[Theorem IV.5.3]{LaSoUr1968}.

The minimum of $w$ on $\overline{\Omega}_T$ is attained at the parabolic boundary. Since the initial data are positive a non-positive minimum can only be attained in a point $(x_0,t_0)\in \Gamma_T$. Then the Robin boundary condition implies that $w(x_0,t_0)\geq v(x_0,t_0)> 0$. This proves the positivity of $w$.

To prove that $\calL$ is continuous we consider $\tilde{v}_1,\tilde{v}_2\in C^{1+\lambda,\frac{1+\lambda}{2}}(\Gamma_T)$ and let $w_j=\calL(\tilde{v}_j),j=1,2$. Since the PDE for $w_j$ is linear in $w_j$ and the boundary condition is linear in $w_j$ and $\tilde{v}_j$ we deduce, from \cite[Theorem IV.5.3]{LaSoUr1968}
\begin{equation*}
\|w_1-w_2\|_{C^{2+\lambda,\frac{2+\lambda}{2}}(\Omega_T)}\leq C\|\tilde{v}_1-\tilde{v}_2\|_{C^{1+\lambda,\frac{1+\lambda}{2}}(\Omega_T)}
\end{equation*}
\hfill$\square$
\end{pf}

\section{Full system}\label{sub:Full}
In this section we will prove the existence, regularity and uniqueness of solutions of the full system.

\begin{thm}\label{thm:ex1-shorttime}
 Let Assumption \ref{ass:Omega} hold. Then there exist $T>0$ and a positive solution $(u,v,w)$ of Problem \ref{sys:Ana} with  $u,v\in C^{2+\lambda,\frac{2+\lambda}{2}}(\Gamma_T)$, $w\in C^{2+\lambda,\frac{2+\lambda}{2}}(\Omega_T)$,
 satisfying the initial conditions \eqref{eq:ini-uv},\eqref{eq:ini-w}.
\end{thm}

\begin{pf}
We will use the Schauder fixed point theorem in the set
\begin{equation*}
  X_T:=\Big\{w\in C^{1+\lambda,\frac{1+\lambda}{2}}(\Omega_T) :\|w\|_{C^{1+\lambda,\frac{1+\lambda}{2}}(\Omega_T)} \leq \Lambda,\, w(\cdot,0)=w_0\text{ in }\Omega\Big\}.
\end{equation*}
Here $T>0$ and $\Lambda>0$ will be chosen below.

Let us consider the solution operators for the bulk and surface systems,
\begin{equation*}
 \calL:C^{1+\lambda,\frac{1+\lambda}{2}}(\Gamma_T) \rightarrow C^{2+\lambda,\frac{2+\lambda}{2}}(\Omega_T),\quad
 \tilde v\mapsto w,
\end{equation*}
where $w$ is the solution of Problem \ref{sys:Bulk}, and
\begin{equation*}
 \calM: C^{\lambda,\frac{\lambda}{2}}(\Omega_T)\rightarrow (C^{2+\lambda,\frac{2+\lambda}{2}}(\Gamma_T))^2,\ \tilde w\mapsto (u,v),
\end{equation*}
where $(u,v)$ denotes the solution of Problem \ref{sys:Surf}.

First of all, $\calM$ is well defined due to Theorem \ref{thm:ex2}, while $\calL$ is well defined due to Theorem \ref{thm:ex1}.

We will obtain the solution to Problem \ref{sys:Ana} by a fixed point of the map
\begin{equation*}
 A:X_T\rightarrow C^{2+\lambda,\frac{2+\lambda}{2}}(\Omega_T),\quad \tilde w\mapsto w:=\calL(v),\text{ where } (u,v)= \calM(\tilde w).
\end{equation*}
Any fixed point $w$ of $A$, together with the pair $(u,v)=\calM(w)$ yields a solution of Problem \ref{sys:Ana}. In order to achieve that $A$ maps $X_T$ into itself we have to choose $T$ sufficiently small and $\Lambda$ sufficiently large.

Let $\tilde w\in X_T$ be arbitrary and set $(u,v)=\calM(\tilde w)$, $w=\calL(v)$.
We need to ensure that
\begin{equation*}
 \|w\|_{C^{1+\lambda,\frac{1+\lambda}{2}}(\Omega_T)}\leq \Lambda.
\end{equation*}
By Theorem \ref{thm:ex2} and Theorem \ref{thm:ex1}
\begin{align*}
 \|w\|_{C^{2+\lambda,\frac{2+\lambda}{2}}(\Omega_T)} &\leq C(\lambda,n,\Omega,T)\Big(1+\Lambda_0+\|v\|_{C^{1+\lambda, \frac{1+\lambda}{2}}(\Gamma_T)}\Big)
 \notag\\
 &\leq  C(\lambda,n,\Omega,T)\Big(1+\Lambda_0+C(\lambda,n,\Omega,T) \big(1+ \Lambda_0 + \|\tilde w\|_{C^{\lambda,\frac{\lambda}{2}}(\Omega_T)}\big)\Big)
 \notag\\
 &\leq C(\lambda,n,\Omega,T)(1+\Lambda_0+\Lambda),
\end{align*}
where the constants $C(\lambda,n,\Omega,T)$ remain bounded with $T\downarrow 0$. We next use \cite[Proposition 16]{AnRoe17} (see also \cite[Lemma B.1]{Angu10}) and deduce that
\begin{equation*}
 \|w\|_{C^{1+\lambda,\frac{1+\lambda}{2}}(\Omega_T)} \leq C(\lambda,n,\Omega,T)T^\delta(1+\Lambda_0+\Lambda) + \Lambda_0, 
\end{equation*}
with $\delta=\delta(\lambda)>0$ and a constant $C(\lambda,n,\Omega,T)$ that remains bounded with $T\downarrow 0$.
Choosing $\Lambda$ sufficiently large relative to $\Lambda_0$ and then $T>0$ sufficiently small yields that $w\in X_T$.

Since $A:X_T\to C^{2+\lambda,\frac{2+\lambda}{2}}(\Omega_T)$ we also have that $A:X_T\to X_T$ compactly.

Since $\calM:C^{\lambda,\frac{\lambda}{2}}(\Omega_T)\to C^{2+\lambda,\frac{2+\lambda}{2}}(\Gamma_T)^2$ is continuous by Theorem \ref{thm:ex2}
and $\calM: C^{1+\gamma, \frac{1+\gamma}{2}}(\Gamma_T)\to C^{2+\lambda,\frac{2+\lambda}{2}}(\Omega_T)$ is linear and bounded by Theorem \ref{thm:ex1} we deduce that $A:X_T\to X_T$ is continuous.

Now by Schauder's fixed point Theorem \cite[Theorem 1.9]{Roub13} there exists $w\in X_T$ with $A(w)=w$ and we obtain with $(u,v)=\calL(w)$ a solution $(u,v,w)$ of Problem \ref{sys:Ana}.

\hfill$\square$
\end{pf}

\begin{prop}[Uniqueness]
\label{prop:unique}
Let Assumption \ref{ass:Omega} hold. Then there exists at most one solution of Problem \ref{sys:Ana}.
\end{prop}

\begin{pf}
Consider two classical solutions $(u,v,w)$ and $(\tilde{u},\tilde{v},\tilde{w})$ of Problem \ref{sys:Ana}.

The PDEs for $u,v,w$ and $\tilde{u},\tilde{v},\tilde{w}$ imply
\begin{equation*}
 \begin{split}
  \frac{\mathsf{d}}{\mathsf{d}t}&\left(\frac{1}{2}\left(\int_\Gamma \left((u-\tilde{u})^2+\tau_s(v-\tilde{v})^2\right)+\int_\Omega \tau_b(w-\tilde{w})^2\right)\right)\\
  =&\int_\Gamma \partial_t(u-\tilde{u})(u-\tilde{u})+\int_\Gamma \tau_s\partial_t(v-\tilde{v})(v-\tilde{v})+\int_\Omega \tau_b\partial_t(w-\tilde{w})(w-\tilde{w})\\
  =&-\varepsilon^2\int_\Gamma|\nabla (u-\tilde{u})|^2-\int_\Gamma(u-\tilde{u})^2+\int_\Gamma(\sigma-\sigma)(u-\tilde{u})+\int_\Gamma\left(\frac{u^p}{v^q}-\frac{\tilde{u}^p}{\tilde{v}^q}\right)(u-\tilde{u})\\
  &-D_s\int_\Gamma|\nabla (v-\tilde{v})|^2-(1+K)\int_\Gamma(v-\tilde{v})^2+K\int_\Gamma(w-\tilde{w})(v-\tilde{v})\\&+\int_\Gamma\left(\frac{u^r}{v^s}-\frac{\tilde{u}^r}{\tilde{v}^s}\right)(v-\tilde{v})
  -D_b\int_\Omega|\nabla(w-\tilde{w})|^2+K\int_\Gamma(v-\tilde{v})(w-\tilde{w})\\&-K\int_\Gamma(w-\tilde{w})^2-\int_\Omega(w-\tilde{w})^2\\
  \leq&C\left(\int_\Gamma\left((u-\tilde{u})^2+(v-\tilde{v})^2\right)+\int_\Omega (w-\tilde{w})^2\right).
 \end{split}
\end{equation*}
Hence $u=\tilde{u},\ v=\tilde{v}$ and $w=\tilde{w}$ follows in view of the Gronwall inequality, since the initial data for the two solutions $(u,v,w)$ and $(\tilde{u},\tilde{v},\tilde{w})$ are equal.

\hfill$\square$
\end{pf}

In the following we consider the solution obtained in Theorem \ref{thm:ex1-shorttime} on the maximal time interval of existence $(0,T)$ and prove that $T=\infty$.
From Theorem \ref{thm:ex1-shorttime} and a continuation argument we already obtain that $u,v\in C^{2+\lambda,\frac{2+\lambda}{2}}_{\loc}(\Gamma_T)$,
$w\in C^{2+\lambda,\frac{2+\lambda}{2}}_{\loc}(\Omega_T)$ hold.

We start again with some lower bounds for the solution.

\begin{lem}
Let Assumption \ref{ass:Omega} hold and let $(u,v,w)$ be a classical solution of Problem \ref{sys:Ana} on $[0,T)$. Then
\begin{align*}
 \inf_{\Gamma_T}u\geq m_u,\quad
 \inf_{\Gamma_T}v\geq m_v
\end{align*}
holds with constants $m_u\geq 0,\, m_v>0$ independent of $T$.
\end{lem}

\begin{pf}
The proof of this Lemma follows immediately from Lemma \ref{lem:Unten2}. Since $(u,v,w)$ is a positive solution of Problem \ref{sys:Ana}, we observe that $(u,v)$ solves Problem \ref{sys:Surf} with $\tilde w=w$. Then Lemma \ref{lem:Unten2} implies that $u$ and $v$ are bounded from below by constants $m_u,m_v$.
\hfill$\square$
\end{pf}

We next prove an a-priori bound in Hölder spaces, which yields the global existence of solution.
\begin{lem}\label{lem:Lk2}
Let Assumption \ref{ass:Omega} hold and let $(u,v,w)$ be a classical solution of Problem \ref{sys:Ana} on $[0,T)$. Let $k\geq 1$ be arbitrary. Then there exists a positive constant $C_k$ independent of $D_s,\, D_b$ and $T$ such that
\begin{align*}
\|u(\cdot,t)\|_{L^k(\Gamma)}+
\|v(\cdot,t)\|_{L^k(\Gamma)}+
\|w(\cdot,t)\|_{L^k(\Omega)}\leq C_k
\end{align*}
holds for all $0<t<T$.
\end{lem}

\begin{pf}
The estimate for $u$ follows by Lemma \ref{lem:Lk1}, with $\tilde w$ replaced by $w$ in Problem \ref{sys:Surf}.

We next prove the corresponding estimates for $v,w$.
The PDE for $v$ implies the equality
\begin{equation*}
\begin{split}
\tau_s \frac{\mathsf{d}}{\mathsf{d}t}\|v\|_{L^k(\Gamma)}^k=&-k\|v\|^k_{L^k(\Gamma)}-k(k-1)D_s\int_\Gamma v^{k-2}|\nabla_\Gamma v|^2\\&-k\int_\Gamma K(v-w)v^{k-1}+k\varepsilon^{-1}\int_\Gamma\frac{u^r}{v^s}v^{k-1},
\end{split}
\end{equation*}
while the PDE and boundary condition for $w$ implies
\begin{equation*}
\tau_b\frac{\mathsf{d}}{\mathsf{d}t}\|w\|_{L^k(\Omega)}^k=-k\|w\|^k_{L^k(\Omega)}-k(k-1)D_s\int_\Omega w^{k-2}|\nabla w|^2+k\int_\Gamma w^{k-1}K(v-w).
\end{equation*}
Since $v^{k-2}|\nabla v|^2$ and $w^{k-2}|\nabla w|^2$ are positive, this leads to the inequality
\begin{equation*}
 \begin{split}
  \tau_s \frac{\mathsf{d}}{\mathsf{d}t}\|v\|_{L^k(\Gamma)}^k+\tau_b\frac{\mathsf{d}}{\mathsf{d}t}\|w\|_{L^k(\Omega)}^k
  \leq
  &-k\left(\|v\|^k_{L^k(\Gamma)}+
  \|w\|^k_{L^k(\Omega)}\right)+\frac{k}{\varepsilon}\int_\Gamma v^{k-1}\frac{u^r}{v^s}\\&-kK\int_\Gamma(v-w)(v^{k-1}-w^{k-1}).
 \end{split}
\end{equation*}
The monotonicity of the mapping $x\mapsto x^{k-1}$ implies $-kK\int_\Gamma(v-w)(v^{k-1}-w^{k-1})>0$, hence the Young inequality implies
\begin{equation*}
 \begin{split}
  \tau_s \frac{\mathsf{d}}{\mathsf{d}t}\|v\|_{L^k(\Gamma)}^k+\tau_b\frac{\mathsf{d}}{\mathsf{d}t}\|w\|_{L^k(\Omega)}^k
  \leq
  & -k\Big(\|v\|^k_{L^k(\Gamma)}
  +\|w\|^k_{L^k(\Omega)}\Big)+\\
  &\frac{k}{\varepsilon}\Big(\delta\int_\Gamma v^k+C_k(\delta)\int_\Gamma\frac{u^{kr}}{v^{ks}}\Big).
 \end{split}
\end{equation*}
With $\delta>0$ chosen sufficiently small we get
\begin{align*}
 &\frac{\mathsf{d}}{\mathsf{d}t}\Big( \tau_s\|v\|_{L^k(\Gamma)}^k+\tau_b\|w\|_{L^k(\Omega)}^k\Big)
 \\
 &\qquad\qquad\qquad
 \leq -\frac{k-1}{\min\{\tau_s,\tau_b\}}\left(\tau_s\|v\|^k_{L^k(\Gamma)}
 +\tau_b\|w\|^k_{L^k(\Omega)}\right)+C(\tau_s,\tau_b)\|u\|_{L^{kr}(\Gamma)}^{kr}.
\end{align*}
From this estimate we conclude as in the end of the proof of Lemma \ref{lem:Lk3} the required bounds for $v,w$.
\hfill$\square$
\end{pf}

\begin{lem}\label{lem:Hoel}
Let Assumption \ref{ass:Omega} hold and let $(u,v,w)$ be a classical solution of Problem \ref{sys:Ana} on $[0,T)$. Then
\begin{equation*}
 u,v\in C^{2+\lambda,\frac{2+\lambda}{2}}(\Gamma_T),\quad
 w\in C^{2+\lambda,\frac{2+\lambda}{2}}(\Omega_T)
\end{equation*}
holds.
\end{lem}

\begin{pf}
We write the PDE for $u$ in the form
\begin{equation*}
 \partial_t u-\varepsilon^2\Delta_\Gamma u=f_1
\end{equation*}
with $f_1=\sigma+\frac{u^p}{v^q}$. Since $u,v$ are in $L^k(0,T;L^k(\Gamma))$ for all $k$ by Lemma \ref{lem:Lk2} and $v$ is greater than a positive number, $f_1$ is within all classes $L^N(0,T;L^N(\Gamma))$.

Regularity theory in parabolic Sobolev spaces \cite[Theorem IV.9.1]{LaSoUr1968} implies $u\in W^{2,1}_N(\Gamma_T)$. Choosing $N>n$ sufficiently large the embedding theorem \cite[Lemma II.3.3]{LaSoUr1968} yields $u\in C^{1+\lambda,\frac{1+\lambda}{2}}(\Gamma_T)$.

With analogous arguments we also derive $v\in C^{1+\lambda,\frac{1+\lambda}{2}}(\Gamma_T)$.
Then regularity theory in Hölder spaces \cite[Theorem IV.5.3]{LaSoUr1968} implies $w\in C^{2+\lambda,\frac{2+\lambda}{2}}(\Omega_T)$.
Applying the corresponding results to $u$ and $v$ yields by \cite[Theorem IV.5.2]{LaSoUr1968} that $u,v\in C^{2+\lambda,\frac{2+\lambda}{2}}(\Gamma_T)$.
\hfill$\square$
\end{pf}

\begin{pot3}
The local existence of solutions has been proved in Theorem \ref{thm:ex1-shorttime}, uniqueness of the maximal solution in Proposition \ref{prop:unique}. Lemma \ref{lem:Hoel} proves an a-priori bound on Hölder norms, which by standard arguments proves the existence of solutions on $(0,\infty)$. We finally prove global-in-time Hölder estimates.

We therefore argue for each of the PDEs in Problem \ref{sys:Ana} separately. Since the arguments are all similar we only consider the PDE for $v$. Let $f_2=-Kw+\eps^{-1}\frac{u^r}{v^s}$, then
\begin{equation*}
 \tau_s\partial_t v -D_s\Delta v = f_2.
\end{equation*}
Define for $k\in\N$ the functions $v_k:\Gamma\times [0,3]\to\R$, $v_k(\cdot,t):= \eta(t)v(\cdot,t+k)$, where $\eta\in C^\infty_c((0,3))$ is a fixed cut-off function with $\eta=1$ on $[1,2]$. We then deduce
\begin{align*}
 \tau_s\partial_t v_k(\cdot,t) -D_s\Delta v_k(\cdot,t) &= \eta(t) f_2(\cdot,t+k) + \eta'(t)v(\cdot,t)
 \quad\text{ in }\Gamma\times (0,3),\\
 v_k(\cdot,0) &= 0 \quad\text{ on }\Gamma.
\end{align*}
Then the right-hand side is bounded in $L^N(\Gamma\times (0,3))$ independently of $k$. Parabolic regularity theory then implies that $v_k$ is bounded in $W^{2,1}_N(\Gamma\times (0,3))$ independently of $k$.
Then $v$ is uniformly bounded in all $W^{2,1}_N(\Gamma\times (k+1,k+2))$, $k\in\N$, hence also in $C^{\lambda,\frac{\lambda}{2}}(\Gamma\times [k+1,k+2])$ by the embedding theorem \cite[Lemma II.3.3]{LaSoUr1968}.
Using this information in the PDE for $v_k$ implies that $v_k$ is bounded in $C^{2+\lambda,1+\frac{\lambda}{2}}_N(\Gamma\times [0,3])$ independently of $k$, which eventually proves that $v\in C^{2+\lambda,1+\frac{\lambda}{2}}(\Gamma\times [0,\infty))$.
\end{pot3}

\section{Discretization of the reduced system}\label{sec:3}
In the limit $D_b\to\infty$ of infinite bulk diffusivity we formally obtain the reduced system that was formulated in Problem \ref{sys:Num}, 
\begin{align*}
 \partial_tu&=\varepsilon^2\Delta_\Gamma u-u+\frac{u^p}{v^q}+\sigma && \text{on}\ \Gamma\times(0,T),\\
 \tau_s\partial_tv&=D_s\Delta_\Gamma v-(1+K)v+\frac{K}{|\Omega|}w+\varepsilon^{-1}\frac{u^r}{v^s} && \text{on}\ \Gamma\times(0,T),\\
 \tau_b\frac{\mathsf{d}}{\mathsf{d}t} w&=\left(1-K\frac{|\Gamma|}{|\Omega|}\right)w+K\int_\Gamma v &&\text{in}\ (0,T).
\end{align*}
We discretize this system using a finite element approximation on a triangulation of a sphere as described in~\cite{Sattar2015}. For a general description of the Finite Element Method we refer the reader to \cite{Kuzmin2011}. The nonlinearities are linearized in a positivity-preserving manner using the Patankar scheme (cf. \cite{Patank1980,BuDeMe2003}). The Patankar linearization of a nonlinear term $h(u)$ is defined by $\tilde{c}(u)u$, where $\tilde{c}(u)=\frac{h(u)}{u}$ is a nonlinear reaction coefficient. This representation yields
\begin{align}
\partial_tu&=\varepsilon^2\Delta_\Gamma u-u+\frac{u^p}{v^{q+1}}v+\sigma && \text{on}\ \Gamma\times(0,T),\label{eq:rsys2a}\\
\tau_s\partial_tv&=D_s\Delta_\Gamma v-(1+K)v+\frac{K}{|\Omega|}w+\varepsilon^{-1}\frac{u^{r-1}}{v^s}u && \text{on}\ \Gamma\times(0,T),\\
\tau_b\frac{\mathsf{d}}{\mathsf{d}t} w&=\left(1-K\frac{|\Gamma|}{|\Omega|}\right)w+K\int_\Gamma v &&\text{in}\ (0,T).\label{eq:rsys2c}
\end{align}
Discretization in space using linear finite elements
leads to the semi-discrete system
\begin{align}
\tau_b\frac{\mathsf{d}}{\mathsf{d} t}w=&- \left(1+K\frac{|\Gamma|}{|\Omega|}\right)w+K\sum_i(G(\vec{v}))_i,\label{eq:rsys3a}\\
\begin{pmatrix}M&&0\\0&&M\end{pmatrix}\begin{pmatrix}\frac{\mathsf{d}}{\mathsf{d} t}\vec{u}\\\frac{\mathsf{d}}{\mathsf{d} t}\vec{v}\end{pmatrix}=&-\begin{pmatrix}A_u&&0\\ 0&&A_v\end{pmatrix}\begin{pmatrix}\vec{u}\\\vec{v}\end{pmatrix}-\begin{pmatrix}B_{u}(\vec{u},\vec{v})&&0\\0&&B_{v}(\vec{u},\vec{v})\end{pmatrix}\begin{pmatrix}\vec{u}\\\vec{v}\end{pmatrix}\nonumber\\&+ \begin{pmatrix}0&&C_{u}(\vec{u},\vec{v})\\C_{v}(\vec{u},\vec{v})&&0\end{pmatrix}\begin{pmatrix}\vec{u}\\\vec{v}\end{pmatrix}+\begin{pmatrix}D\\E(w)\end{pmatrix}.\label{eq:rsys3c}
\end{align}
For discretization in time we will use two positivity-preserving methods. The first method will lead to a first order accurate discretization, while the second one will produce a flux-corrected second order discretization.

\subsection{First order discretization}
Using a first order operator splitting ansatz, we will solve the ODE and the semi-discrete problem for the coupled reaction-diffusion equations on the surface sequentially at each time step. For the discretization of the ODE we use the implicit Euler scheme. The fully discrete counterpart of the ODE step reads
\begin{equation}
\tau_b\frac{w^{n+1}-w^n}{\Delta t}=-\left(1+K\frac{|\Gamma|}{|\Omega|}\right)w^{n+1}+K\sum_i(G(\vec{v}^n))_i
\end{equation}
with initial condition $w^0=w_0$. In this case the implicit Euler scheme is positivity preserving. If $\vec{v}^n$ is positive, then, after solving for $w^{n+1}$, all coefficients are positive. For discretization in time of the surface sub-system we choose the positivity preserving Patankar-Euler scheme (cf. \cite{BuDeMe2003}). This discretization leads to the system
\begin{equation}\begin{split}
\begin{pmatrix}M+\Delta t(A_u+B_u(\vec{u}^n,\vec{v}^n))&&-\Delta tC_u(\vec{u}^n,\vec{v}^n)\\-\Delta tC_v(\vec{u}^n,\vec{v}^n)&&M+\Delta t(A_v+B_v(\vec{u}^n,\vec{v}^n))\end{pmatrix}\begin{pmatrix}\vec{u}^{n+1}\\\vec{v}^{n+1}\end{pmatrix}\\=\begin{pmatrix}M&&0\\0&&M\end{pmatrix}\begin{pmatrix}\vec{u}^n\\\vec{v}^n\end{pmatrix}+ \begin{pmatrix}D\\E(w^{n+1})\end{pmatrix}\end{split}\label{step2}
\end{equation}
with the initial condition $(\vec{u}^0,\vec{v}^0)^T=(\vec{u}_0,\vec{v}_0)^T$. The lumped counterpart of the mass matrix $M$ is defined by $ \tilde{M}=\{\tilde{m}_i\Delta_{ij}\}$, where $\tilde{m}_{i}=\sum_j m_{ij}$. The lumped reactive matrices $\tilde{B}_u,\tilde{B}_v,\tilde{C}_u,\tilde{C}_v$ are defined similarly. The lumped version of system \eqref{step2} is positivity preserving provided that the off-diagonal entries of $A_u$ and $A_v$ are non-positive. For piecewise-linear finite element discretizations of the Laplace operator, this requirement is satisfied under certain assumptions regarding the triangulation. In view of these considerations, we consider the lumped Euler-Patankar scheme
\begin{align*}
	w^{n+1}=\frac{\tau_b}{\tau_b+\Delta t\left(1+\frac{(K|\Gamma|)}{|\Omega|}\right)}w^n+\frac{\Delta tK}{\tau_b+\Delta t\left(1+\frac{(K|\Gamma|)}{|\Omega|}\right)}\sum_i(G(\vec{v}^n))_i,\\
	\begin{pmatrix}\tilde{M}+\Delta t(A_u+\tilde{B}_u(\vec{u}^n,\vec{v}^n))&&-\Delta t\tilde{C}_u(\vec{u}^n,\vec{v}^n)\\-\Delta t\tilde{C}_v(\vec{u}^n,\vec{v}^n)&&\tilde{M}+\Delta t(A_v+\tilde{B}_v(\vec{u}^n,\vec{v}^n))\end{pmatrix}\begin{pmatrix}\vec{u}^{n+1}\\\vec{v}^{n+1}\end{pmatrix}\nonumber\\= \begin{pmatrix}\tilde{M}&&0\\0&&\tilde{M}\end{pmatrix}\begin{pmatrix}\vec{u}^n\\\vec{v}^n\end{pmatrix}+ \Delta t\begin{pmatrix}D\\E(w^{n+1})\end{pmatrix}
\end{align*}
which is first-order accurate in time and positivity preserving if the discretization of the Laplace-Beltrami operator produces matrices $A_u$ and $A_v$ with desired properties.

\subsection{Second order discretization}
Let us now construct a time stepping scheme of second order.  We use an operator splitting ansatz again, so we can solve the ODE and the sub-system on the surface separately. To achieve second order accuracy, we use the Strang splitting (cf. \cite{Strang1968,HunVer2003}) which updates the solutions of the two subproblems in the following sequential manner:
\begin{equation}
(u^n, v^n, w^n)\mapsto(u^n, v^n, w^{n+\frac{1}{2}})\mapsto(u^{n+1},v^{n+1},w^{n+\frac{1}{2}})\mapsto(u^{n+1},v^{n+1},w^{n+1}).
\end{equation}
That is, we begin with a half time step for the ODE, followed by a full time step for the surface sub-system and finally a second half time step for the ODE. As the data depending on other variables we always use the most recently computed values. For the discretization of the ODE we use the Patankar-$\theta$-scheme \cite{Gruell2018}. This scheme represents a modification of the second order SSP-Runge-Kutta method
\begin{align}
w^{(1)}&=w^n+\Delta t(a(w^n)-b(w^n)w^n),\\w^{(2)}&=w^{(1)}+\Delta t(a(w^{(1)})-b(w^{(1)})w^{(1)}),\\w^{n+1}&=\frac{1}{2}(w^n+w^{(2)})
\end{align}
for an ODE of the form
\begin{equation}
\frac{\mathsf{d}}{\mathsf{d}t}w=a(w)-b(w)w.
\end{equation}
This scheme is modified by using a correction procedure which has no effect if the scheme is already positivity preserving. The modified scheme reads
\begin{align}
w^{(1)}&=w^n+\Delta t(a(w^n)-b(w^n)((1-\theta_1)w^n+\theta_1w^{(1)})),\\
w^{(2)}&=w^{(1)}+\Delta t(a(w^{(1)})-b(w^{(1)})((1-\theta_2)w^{(1)}+\theta_2w^{(2)})),\\
w^{n+1}&=\frac{1}{2}(w^n+w^{(2)}),
\end{align}
where $\theta_i$, $i=1,2$ should be chosen to ensure that
\begin{equation}
w^{(i)}=w^{(i-1)}+\Delta t (a(w^{(i-1)})-b(w^{(i-1)})((1-\theta_i)w^{(i-1)}+\theta_iw^{(i-1)}))\geq 0
\end{equation}
holds with $w^{(0)}=w^n$. This criterion is satisfied for
\begin{equation}
\theta_i=\max\left\{0,1-\frac{1}{\Delta tb(w^{(i-1)})}-\frac{a(w^{(i-1)})}{b(w^{(i-1)})w^{(i-1)}}\right\}.
\end{equation}
In our case we have $a(w)=\frac{K}{\tau_b}\sum_i(G(\vec{v}))_i$ and $b(w)=\frac{1}{\tau_b}\left(1+K\frac{|\Gamma|}{|\Omega|}\right)$. For the time discretization of the subsystem on the surface we choose the Crank-Nicolson scheme. Since it is generally not positivity preserving, we combine it with the positivity preserving low order time stepping scheme using  Zalesak's flux-corrected transport (FCT) algorithm (cf. \cite{Zalesak1979,Kuzmin2011}). To reduce the number of systems to be solved, a modification of the Crank-Nicolson scheme is performed before the FCT limiter is applied. At the first step, we compute the low order solution $(\vec{u}^L,\vec{v}^L)$ with the Patankar-Euler scheme as described above. Next, the high order solution is computed by solving the system
\begin{equation}
\begin{split}
\begin{pmatrix}M+\frac{\Delta t}{2}A_u&&0\\0&&M+\frac{\Delta t}{2}A_v\end{pmatrix}\begin{pmatrix}\vec{u}^{H}\\\vec{v}^{H}\end{pmatrix}=\begin{pmatrix}M-\frac{\Delta t}{2}A_u&&0\\0&&M-\frac{\Delta t}{2}A_v\end{pmatrix}\begin{pmatrix}\vec{u}^{n}\\\vec{v}^{n}\end{pmatrix}\\+\Delta t\begin{pmatrix}-B_u(\vec{u}^n,\vec{v}^n)&&C_u(\vec{u}^n,\vec{v}^n)\\C_v(\vec{u}^n,\vec{v}^n)&&-B_v(\vec{u}^n,\vec{v}^n)\end{pmatrix}\begin{pmatrix}\vec{u}^L\\\vec{v}^L\end{pmatrix}+\Delta t\begin{pmatrix}D\\ E(w^{n+1})\end{pmatrix}
\end{split}.
\end{equation}
To limit $\vec u^H$ in a way which makes it positivity preserving without losing the discrete conservation property, we decompose the difference between
\begin{equation}
m_i u^L_i=m_i u^n_i + \Delta t(-b_i(\vec{u}^n,\vec{v}^n)u_i^L+c_i(\vec{u}^n,\vec{v}^n)v_i^L)-\Delta t\sum_{j\neq i}a_{ij}(u^L_j-u^L_i)+d_i
\label{eq:PFCT1}
\end{equation}
and
\begin{equation}
\begin{split}
m_i u^H_i=&m_i u^n_i + \Delta t(-b_i(\vec{u}^n,\vec{v}^n)u_i^L+c_i(\vec{u}^n,\vec{v}^n)v_i^L)-\frac{\Delta t}{2}\sum_{j\neq i}a_{ij}(u^H_j-u^H_i)\\&-\frac{\Delta t}{2}\sum_{j\neq i}a_{ij}(u^n_j-u^n_i)-\Delta t\sum_{j\neq i}b_{ij}(\vec{u}^n,\vec{v}^n)(u^L_j-u_i^L)\\&+\Delta t\sum_{j\neq i}c_{ij}(\vec{u}^n,\vec{v}^n)(u^L_j-u_i^L)+\sum_{j\neq i}m_{ij}(u^H_i-u^H_j)\\&-\sum_{j\neq i}m_{ij}(u^n_i-u^n_j)+d_i\end{split}\label{eq:PFCT2}
\end{equation}
into the antidiffusive fluxes
\begin{equation}
\begin{split}
f_{ij}=&a_{ij}\left((u^L_j-u^ L_i)-\frac{1}{2}(u^H_j-u^H_i)-\frac{1}{2}(u^n_j-u^n_i)\right)-b_{ij}(\vec{u}^n,\vec{v}^n)(u^L_j-u_i^L)\\&+c_{ij}(\vec{u}^n,\vec{v}^n)(u^L_j-u_i^L)-\frac{m_{ij}}{\Delta t}\left((u^H_j-u^H_i)-(u^n_j-u^n_i)\right)\end{split}
\end{equation}
such that
\begin{equation}
m_iu^H_i=m_iu^L_i+\Delta t\sum_{j\neq i}f_{ij}.
\end{equation}
To ensure that $u^H_i$ will be bounded by the local maxima and minima of $\vec u^L$,
the fluxes $f_{ij}$ are multiplied by correction factors $\alpha_{ij}\in[0,1]$. We calculate them using Zalesak's algorithm, as presented in \cite{Kuzmin2011}. The flux-corrected solution is given by
\begin{equation}
m_iu^{n+1}_i=m_iu^L_i+\Delta t\sum_{j\neq i}\alpha_{ij}f_{ij}.
\end{equation}
The FCT-constrained approximation $v^{n+1}_i$ is obtained similarly using Zalesak's limiter to enforce positivity preservation. Since we limit the high order solution $\vec v^H$  to preserve the local range of the low order solution $\vec v^L$ computed with a positivity preserving scheme, the flux-corrected solution $(\vec{u}^{n+1},\vec{v}^{n+1})$ is positive as well.

\section{Numerical examples}\label{sec:4}
In what follows, we will investigate the quality of the presented numerical schemes. The domain $\Omega$ that we consider in this numerical study is the three dimensional unit ball. The $2D$ case  was simulated and investigated theoretically by Gomez et al in \cite{GoWaWe2018}. We will compare the results of the $3D$ simulation to expectations that can be inferred from the results of the $2D$ analysis and simulations. In particular, we expect that there is a Turing instability and we can observe pattern formation. The $2D$ simulation results also suggest that formation of one or two spikes may be expected for certain values of the data and of the coefficients. For testing the implementation of our numerical scheme we used the data listed in table \ref{tab:1}. The employed values of $K$ are given by $0.002$ and $200,000.0$. As initial conditions for $w$ and $v$, we used $w_0=0.01$ and $v_0\equiv 0.1$, respectively. The choice of initial data for $u$ was varied in our numerical experiments.
\begin{table}
\centering
\caption{Parameters}
\label{tab:1}
\begin{tabular}{cccc|cccc}
$p$ & $q$ & $r$ & $s$ & $\epsilon$ & $D_s$ & $\tau_s$ & $\tau_b$\\
\hline
$2$ & $4$ & $3$ & $4$ & $0.1$ & $10.0$ & $0.6$ & $0.1$
\end{tabular}
\end{table}

\begin{figure}
\centering
\includegraphics[width=0.3\textwidth]{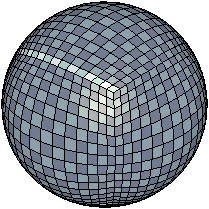}
\caption{Cubed sphere grid which is refined five times to produce the mesh for this computational study}
\label{abb:gra1}
\end{figure}

\begin{figure}
\centering
\includegraphics[width=0.3\textwidth]{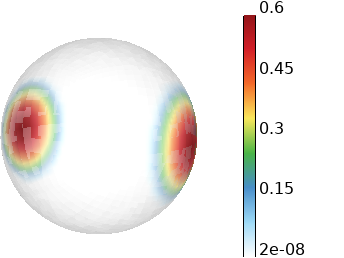}\hspace{0.2\textwidth}
\includegraphics[width=0.3\textwidth]{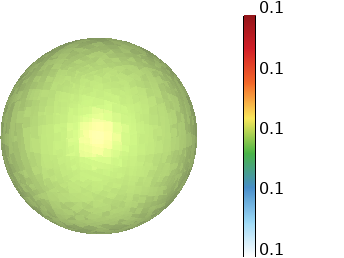}\\
\includegraphics[width=0.3\textwidth]{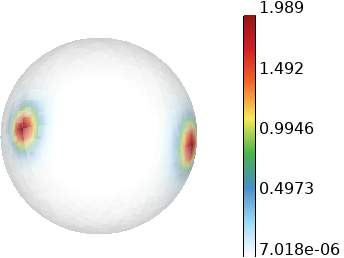}\hspace{0.2\textwidth}
\includegraphics[width=0.3\textwidth]{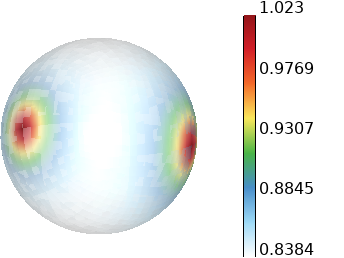}\\
\includegraphics[width=0.3\textwidth]{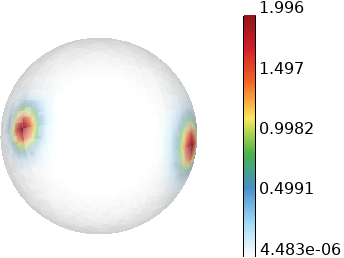}\hspace{0.2\textwidth}
\includegraphics[width=0.3\textwidth]{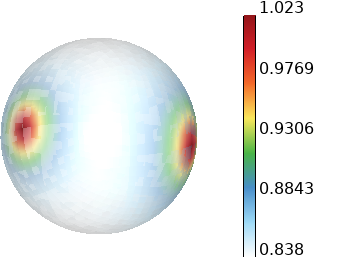}\\
\includegraphics[width=0.3\textwidth]{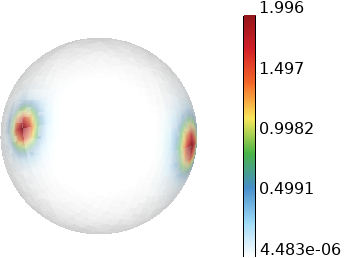}\hspace{0.2\textwidth}
\includegraphics[width=0.3\textwidth]{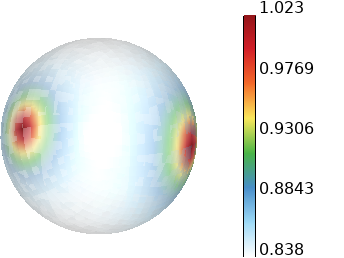}\\
\includegraphics[width=0.3\textwidth]{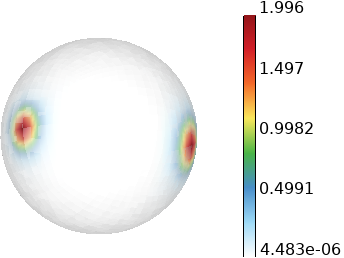}\hspace{0.2\textwidth}
\includegraphics[width=0.3\textwidth]{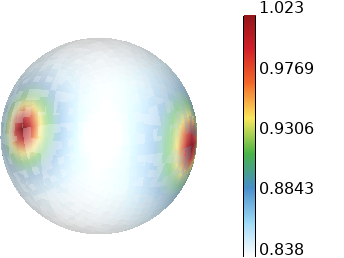}
\caption[Plots of $(u,v)$ with a 180 degree 2-spike initial condition and $K=0.002$]{Plots of $u$ (left) and $v$ (right) for a 180 degree 2-spike initial condition and $K=0.002$. The snapshots correspond to the time instants $t=0,10,20,70,500$ (from top to bottom).}
\label{abb:ZGku}
\end{figure}

\begin{figure}
\centering
\includegraphics[width=0.3\textwidth]{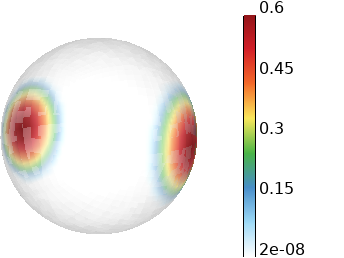}\hspace{0.2\textwidth}
\includegraphics[width=0.3\textwidth]{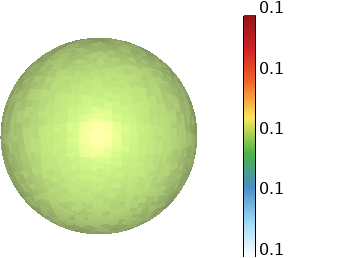}
\includegraphics[width=0.3\textwidth]{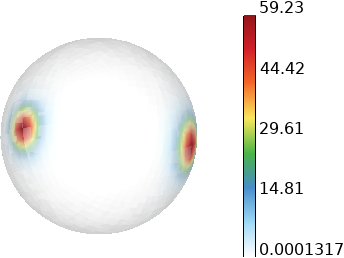}\hspace{0.2\textwidth}
\includegraphics[width=0.3\textwidth]{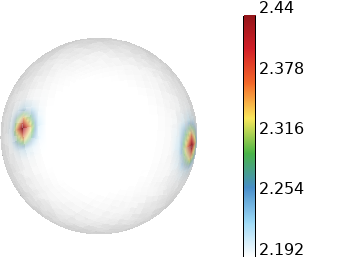}
\includegraphics[width=0.3\textwidth]{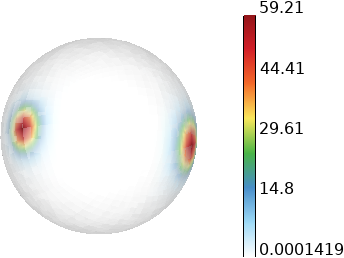}\hspace{0.2\textwidth}
\includegraphics[width=0.3\textwidth]{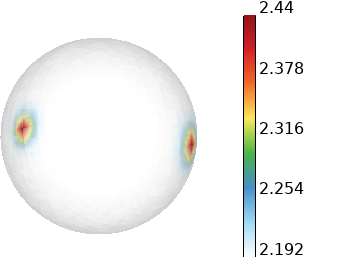}
\includegraphics[width=0.3\textwidth]{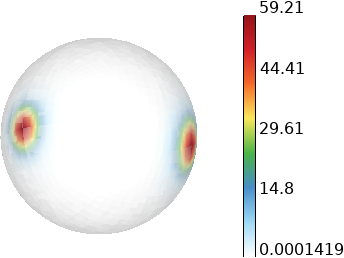}\hspace{0.2\textwidth}
\includegraphics[width=0.3\textwidth]{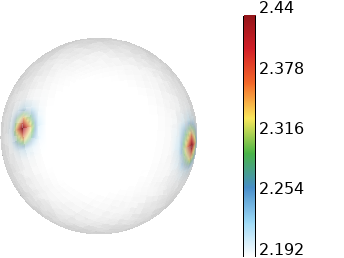}
\includegraphics[width=0.3\textwidth]{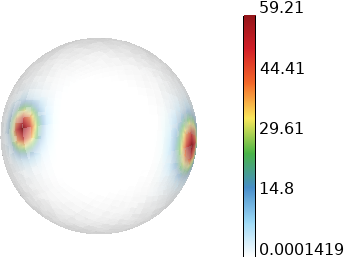}\hspace{0.2\textwidth}
\includegraphics[width=0.3\textwidth]{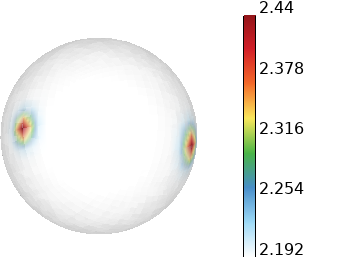}
\caption[Plots of $(u,v)$ with a 180 degree 2-spike initial condition and $K=200000.0$]{Plots of $u$ (left) and $v$ (right) for a 180 degree 2-spike initial condition and $K=200000.0$. The snapshots correspond to the time instants $t=0,10,20,70,500$ (from top to bottom).}
\label{abb:ZGgu}
\end{figure}

\begin{figure}
\centering
\includegraphics[width=0.3\textwidth]{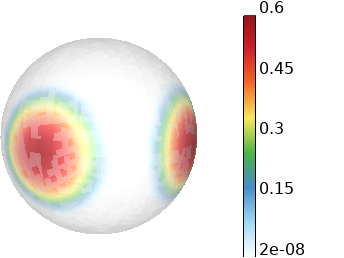}\hspace{0.2\textwidth}
\includegraphics[width=0.3\textwidth]{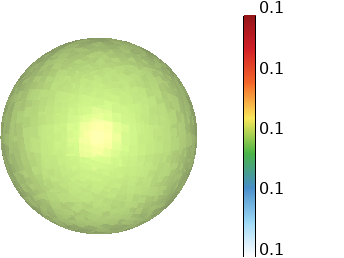}
\includegraphics[width=0.3\textwidth]{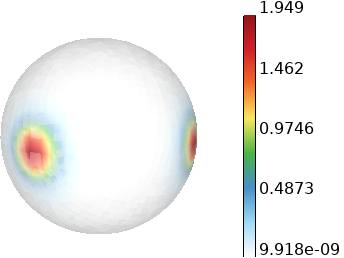}\hspace{0.2\textwidth}
\includegraphics[width=0.3\textwidth]{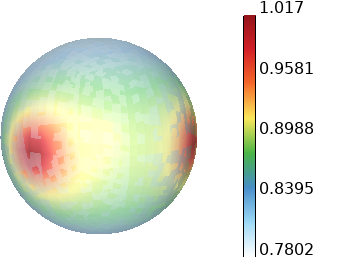}
\includegraphics[width=0.3\textwidth]{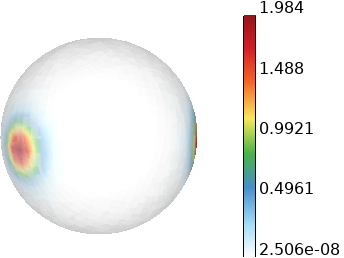}\hspace{0.2\textwidth}
\includegraphics[width=0.3\textwidth]{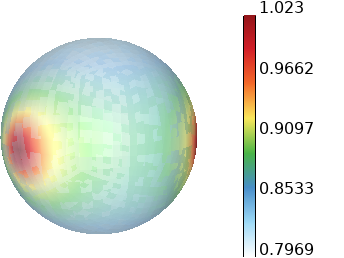}
\includegraphics[width=0.3\textwidth]{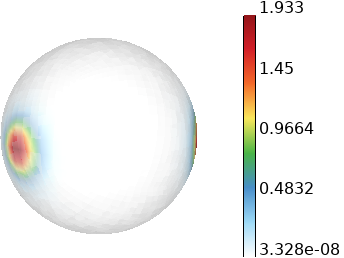}\hspace{0.2\textwidth}
\includegraphics[width=0.3\textwidth]{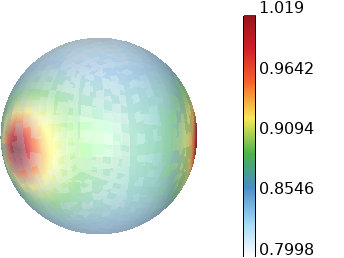}
\includegraphics[width=0.3\textwidth]{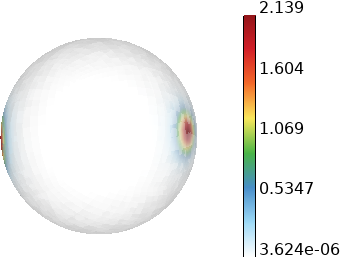}\hspace{0.2\textwidth}
\includegraphics[width=0.3\textwidth]{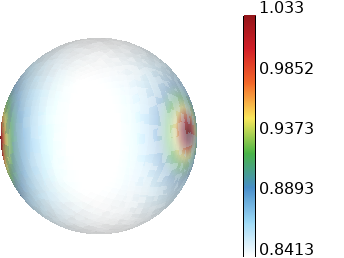}
\caption[Plots of $(u,v)$ with a 90 degree 2-spike initial condition and $K=0.002$]{Plots of $u$ (left) and $v$ (right) for a 90 degree -spike initial condition and $K=0.002$. The snapshots correspond to the time instants $t=0,10,20,70,1000$ (from top to bottom).}
\label{abb:ZEku}
\end{figure}

\begin{figure}
\centering
\includegraphics[width=0.3\textwidth]{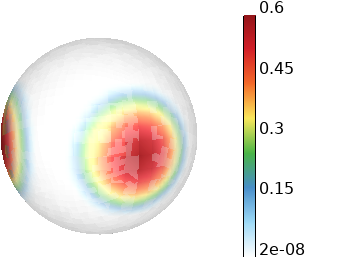}\hspace{0.2\textwidth}
\includegraphics[width=0.3\textwidth]{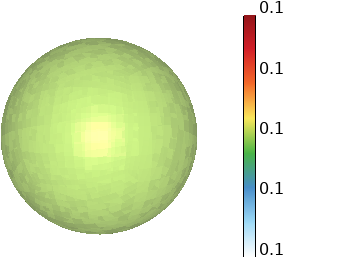}
\includegraphics[width=0.3\textwidth]{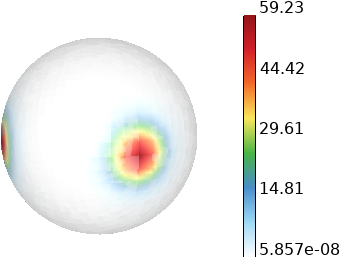}\hspace{0.2\textwidth}
\includegraphics[width=0.3\textwidth]{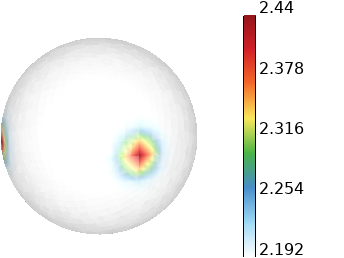}
\includegraphics[width=0.3\textwidth]{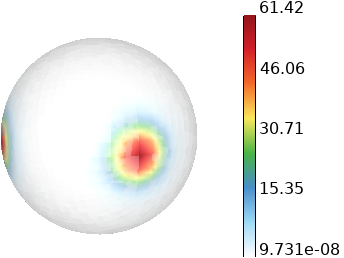}\hspace{0.2\textwidth}
\includegraphics[width=0.3\textwidth]{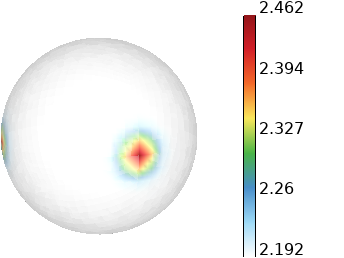}
\includegraphics[width=0.3\textwidth]{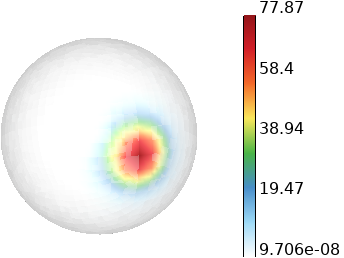}\hspace{0.2\textwidth}
\includegraphics[width=0.3\textwidth]{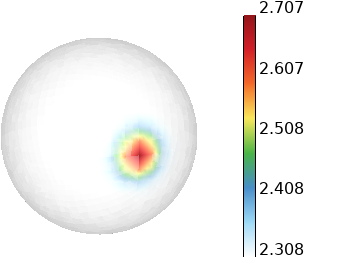}
\includegraphics[width=0.3\textwidth]{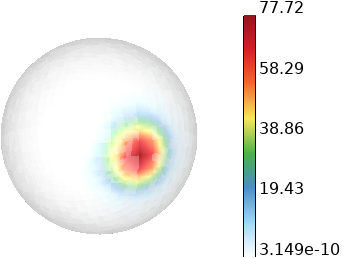}\hspace{0.2\textwidth}
\includegraphics[width=0.3\textwidth]{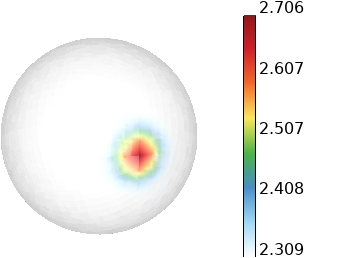}
\caption[Plots of $(u,v)$ with a 90 degree 2-spike initial condition and $K=200000.0$]{Plots of $u$ (left) and $v$ (right) for a 90 degree 6-spike initial condition and $K=200000.0$. The snapshots correspond to the time instants $t=0,10,20,70,500$ (from top to bottom).}
\label{abb:ZEgu}
\end{figure}

\begin{figure}
\centering
\includegraphics[width=0.3\textwidth]{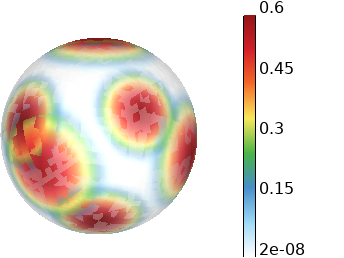}\hspace{0.2\textwidth}
\includegraphics[width=0.3\textwidth]{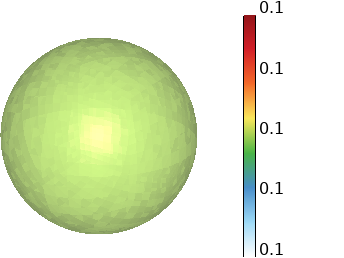}
\includegraphics[width=0.3\textwidth]{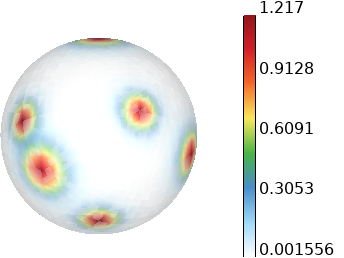}\hspace{0.2\textwidth}
\includegraphics[width=0.3\textwidth]{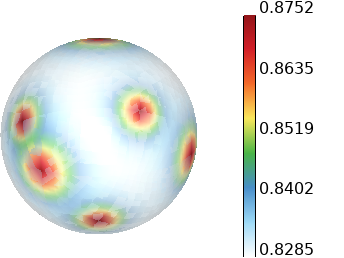}
\includegraphics[width=0.3\textwidth]{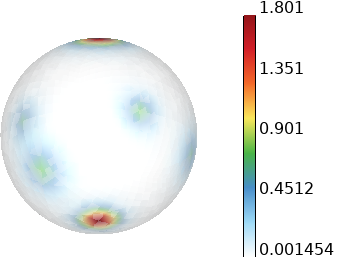}\hspace{0.2\textwidth}
\includegraphics[width=0.3\textwidth]{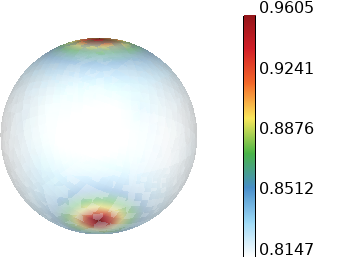}
\includegraphics[width=0.3\textwidth]{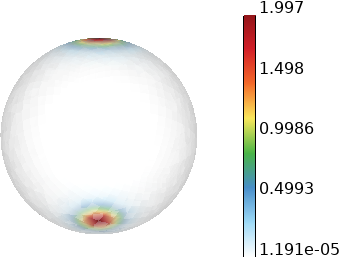}\hspace{0.2\textwidth}
\includegraphics[width=0.3\textwidth]{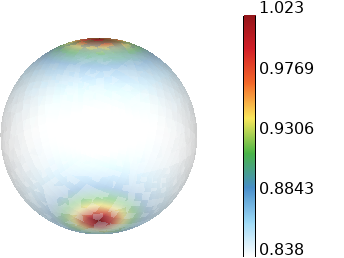}
\includegraphics[width=0.3\textwidth]{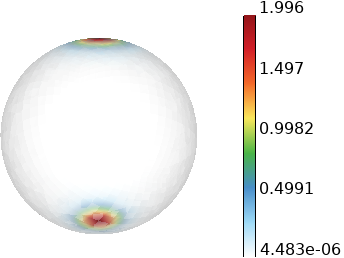}\hspace{0.2\textwidth}
\includegraphics[width=0.3\textwidth]{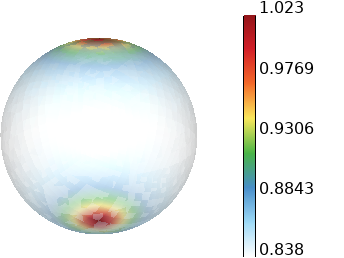}
\caption[Plots of $(u,v)$ with a 90 degree 6-spike initial condition and $K=0.002$]{Plots of $u$ (left) and $v$ (right) for a 90 degree 6-spike initial condition and $K=0.002$. The snapshots correspond to the time instants $t=0,10,20,70,500$ (from top to bottom).}
\label{abb:Sku}
\end{figure}

\begin{figure}
\centering
\includegraphics[width=0.3\textwidth]{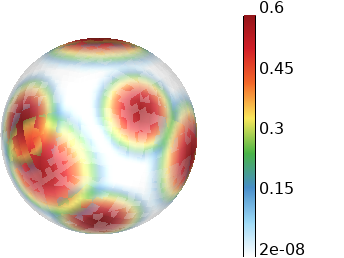}\hspace{0.2\textwidth}
\includegraphics[width=0.3\textwidth]{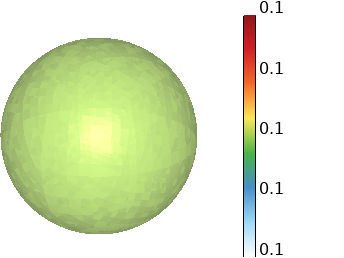}
\includegraphics[width=0.3\textwidth]{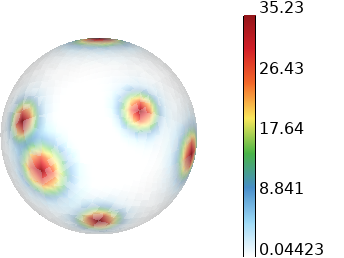}\hspace{0.2\textwidth}
\includegraphics[width=0.3\textwidth]{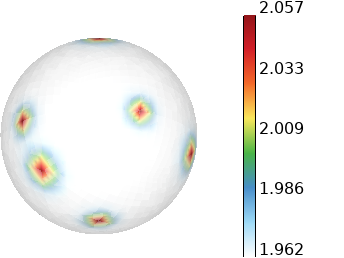}
\includegraphics[width=0.3\textwidth]{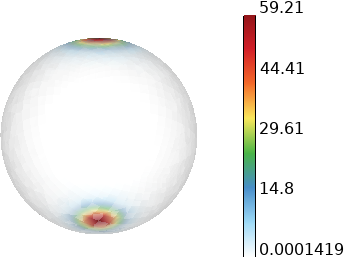}\hspace{0.2\textwidth}
\includegraphics[width=0.3\textwidth]{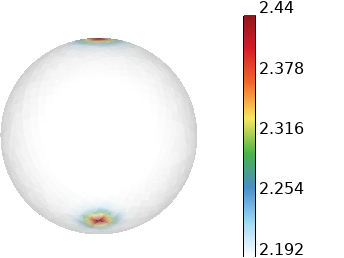}
\includegraphics[width=0.3\textwidth]{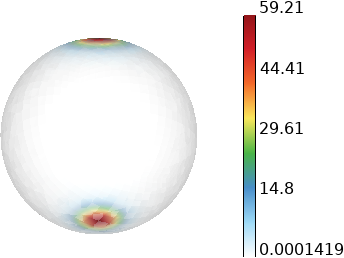}\hspace{0.2\textwidth}
\includegraphics[width=0.3\textwidth]{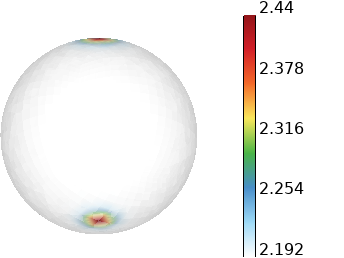}
\includegraphics[width=0.3\textwidth]{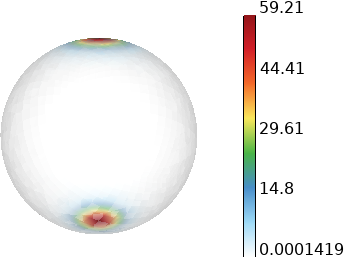}\hspace{0.2\textwidth}
\includegraphics[width=0.3\textwidth]{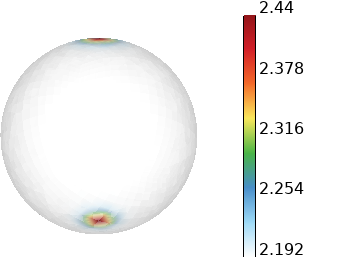}
\caption[Plots of $(u,v)$ with a 90 degree 6-spike initial condition and $K=200000.0$]{Plots of $u$ (left) and $v$ (right) for a 90 degree 6-spike initial condition and $K=200000.0$. The snapshots correspond to the time instants $t=0,10,20,70,500$ (from top to bottom).}
\label{abb:Sgu}
\end{figure}

\begin{figure}
\centering
\includegraphics[width=0.3\textwidth]{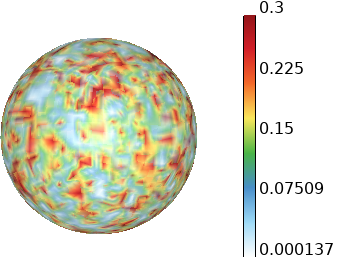}\hspace{0.2\textwidth}
\includegraphics[width=0.3\textwidth]{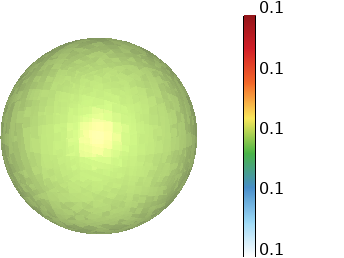}
\includegraphics[width=0.3\textwidth]{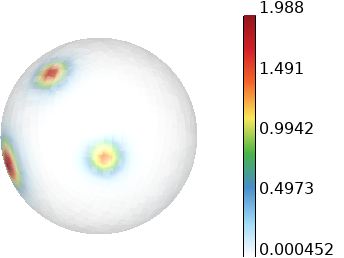}\hspace{0.2\textwidth}
\includegraphics[width=0.3\textwidth]{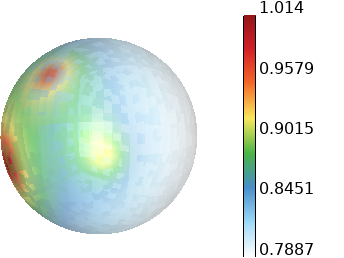}
\includegraphics[width=0.3\textwidth]{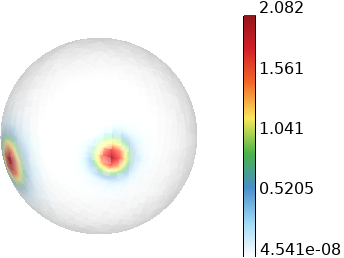}\hspace{0.2\textwidth}
\includegraphics[width=0.3\textwidth]{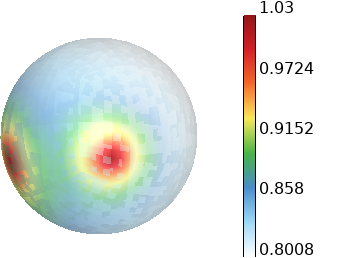}
\includegraphics[width=0.3\textwidth]{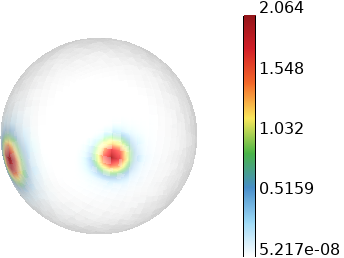}\hspace{0.2\textwidth}
\includegraphics[width=0.3\textwidth]{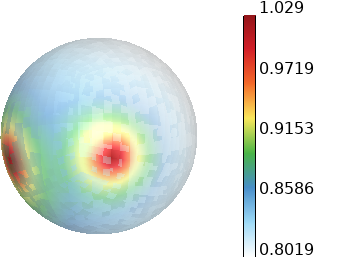}
\includegraphics[width=0.3\textwidth]{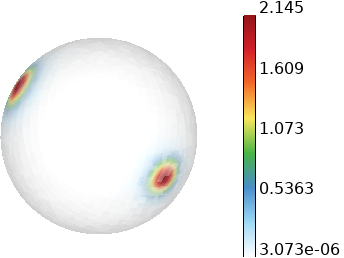}\hspace{0.2\textwidth}
\includegraphics[width=0.3\textwidth]{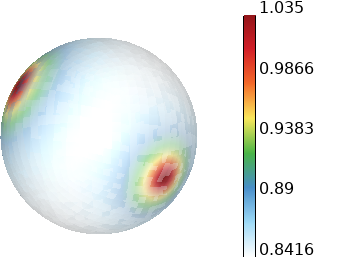}
\caption[Plots of $(u,v)$ with a random initial condition and $K=0.002$]{Plots of $u$ (left) and $v$ (right) for a random initial condition and $K=0.002$. The snapshots correspond to the time instants $t=0,10,20,70,1000$ (from top to bottom).}
\label{abb:Cku}
\end{figure}

\begin{figure}
\centering
\includegraphics[width=0.3\textwidth]{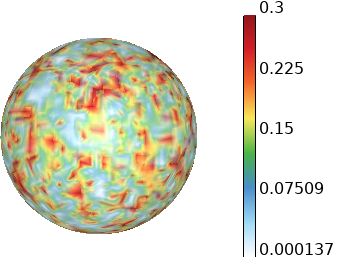}\hspace{0.2\textwidth}
\includegraphics[width=0.3\textwidth]{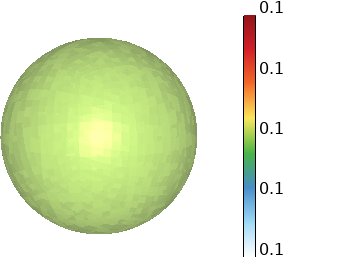}
\includegraphics[width=0.3\textwidth]{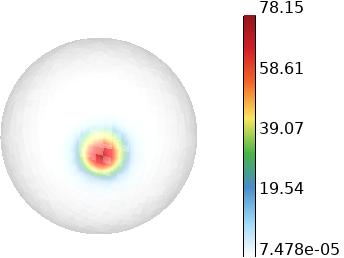}\hspace{0.2\textwidth}
\includegraphics[width=0.3\textwidth]{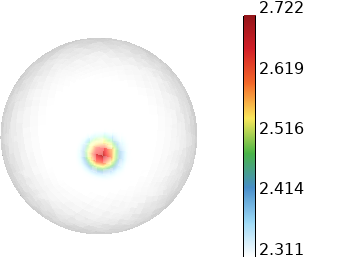}
\includegraphics[width=0.3\textwidth]{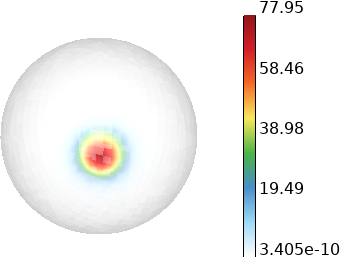}\hspace{0.2\textwidth}
\includegraphics[width=0.3\textwidth]{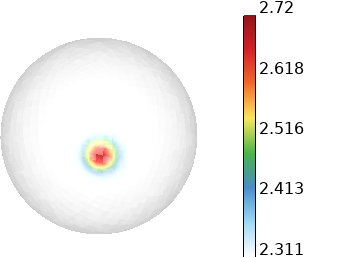}
\includegraphics[width=0.3\textwidth]{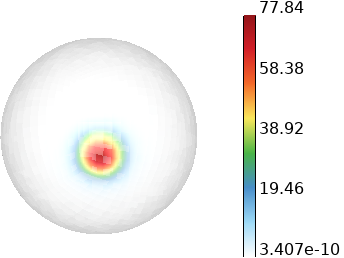}\hspace{0.2\textwidth}
\includegraphics[width=0.3\textwidth]{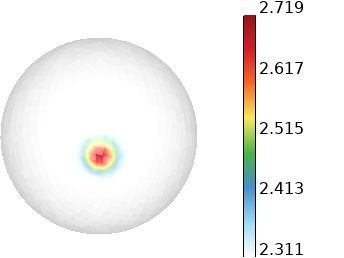}
\includegraphics[width=0.3\textwidth]{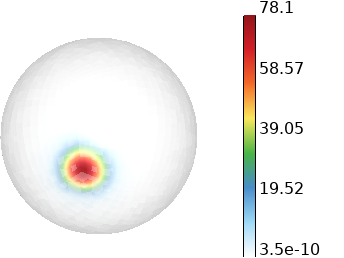}\hspace{0.2\textwidth}
\includegraphics[width=0.3\textwidth]{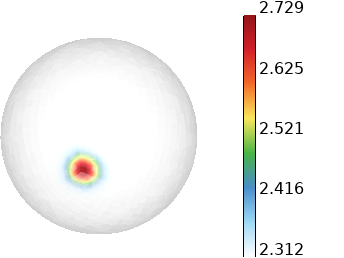}
\caption[Plots of $(u,v)$ with a random initial condition and $K=200000.0$]{Plots of $u$ (left) and $v$ (right) for a random initial condition and $K=200000.0$. The snapshots correspond to the time instants $t=0,10,20,70,1000$ (from top to bottom).}
\label{abb:Cgu}
\end{figure}

The simulations were performed using the time step $\Delta t=10^{-5}$ on the fifth level of uniform refinement for the cubed sphere grid shown in Fig. \ref{abb:gra1}. For the simulations the mfem library (cf. \cite{mfem,mfem-web}) is used, while the visualization is done with GlVis (cf. \cite{glvis-tool}). We have simulated a total of eight different cases corresponding to four different initial conditions $u_0$ and two different values of $K$. The values $K=200,000.00$ and $K=0.002$ were used for each initial condition. The initial shapes of $u$ included different spike patterns and a random distribution. In Table \ref{tab:2} it is shown which initial condition leads to which stationary solution pattern depending on $K$. A symmetric 2-spike initial condition results in a pattern formation that generates two spikes on opposite sides of the sphere. A nonsymmetric 2-spike initial condition produces two spikes which form a 90 degree angle, while the 6-spike initial condition produces one spike on each face of the cube which is projected on the sphere. The minimal and maximal values of $u$ listed in table \ref{tab:2} are calculated at the time $t=500.00$ respectively $t=1000.00$ for the cases marked with (*). The simulation of these eight cases reveals two different patterns in stationary solutions. One of these patterns preserves a single spike, while the other pattern preserves two spikes on opposite sides of the sphere. In this study we always obtained a symmetric 2-spike pattern with $K=0.002$, while the patterns obtained with $K=200,000.00$ were found to depend on the initial data. In figures \ref{abb:ZGku}, \ref{abb:ZEku}, \ref{abb:Sku} and \ref{abb:Cku} we observe the formation of a two spike pattern with minima and maxima of comparable magnitude. In figures \ref{abb:ZGgu} and \ref{abb:Sgu} we also observe the formation of a two spike pattern but the magnitude of the minima and maxima changes. Finally in figures \ref{abb:ZEgu} and \ref{abb:Cgu} we observe the formation of a one spike pattern. These observations imply that the pattern depends on $K$ and on the structure of the initial condition.
\begin{table}
\centering
\caption{Solution patterns at the end of the simulations. In the cases marked with (*) the solutions are not yet stationary, but appear to approach the respective stationary states.}
\label{tab:2}
\begin{tabular}{c|c|c|c|c}
initial condition&$K$&stationary solution& $\min$&$\max$\\
\hline
2-spike symmetric&0.002&2-spike symmetric& $4.483\cdot 10^{-6}$&$1.996$\\
2-spike symmetric&200,000.00&2-spike symmetric& $1.419\cdot 10^{-4}$&$59.21$\\
\hline
2-spike nonsymmetric(*)&0.002&2-spike symmetric& $3.624\cdot 10^{-6}$&$2.139$\\
2-spike nonsymmetric&200,000.00&1-spike& $3.149\cdot 10^{-10}$&$77.72$\\
\hline
6-spike&0.002&2-spike symmetric& $4.483\cdot 10^{-6}$&$1.997$\\
6-spike&200,000.00&2-spike symmetric & $1.419\cdot 10^{-4}$&$59.21$\\
\hline
random(*)&0.002&2-spike symmetric & $3.073\cdot 10^{-6}$&$2.145$\\
random(*)&200,000.00&1-spike& $3.5\cdot 10^{-10}$&$78.1$\\
\end{tabular}
\end{table}

Therefore the results of our simulations are consistent with what we expected based on the published results for the two dimensional case. Furthermore we see that no blow-ups occur in finite time. This is also to be expected in view of the existence of a solution to the full system which we proved in \ref{sub:Full}.

\begin{figure}
\centering
\includegraphics[width=0.3\textwidth]{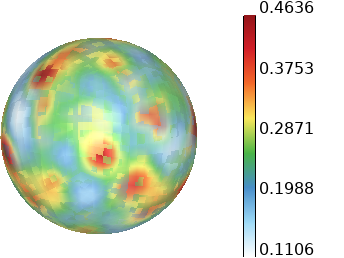}\hspace{0.2\textwidth}
\includegraphics[width=0.3\textwidth]{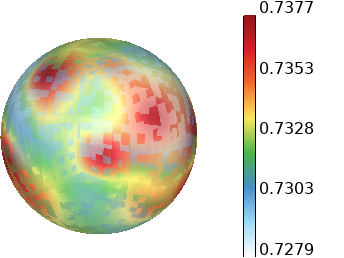}
\includegraphics[width=0.3\textwidth]{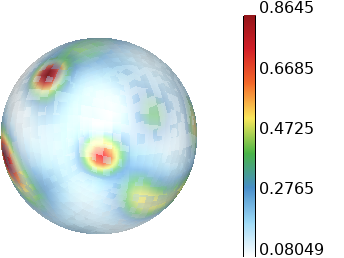}\hspace{0.2\textwidth}
\includegraphics[width=0.3\textwidth]{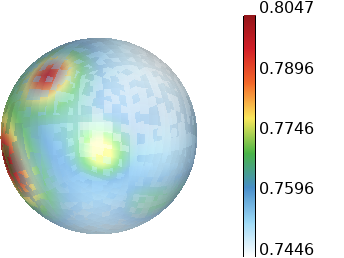}
\caption[Additional plots of $(u,v)$ with a random initial condition and $K=0.002$]{Additional plots of $u$ (left) and $v$ (right) for a random initial condition and $K=0.002$. The snapshots correspond to the time instants $t=1, 3$ (from top to bottom).}
\label{abb:Cku2}
\end{figure}

\begin{figure}
\centering
\includegraphics[width=0.3\textwidth]{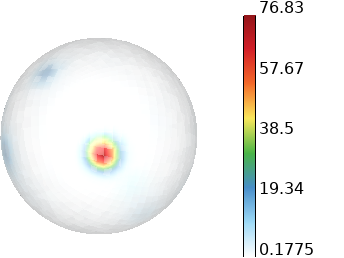}\hspace{0.2\textwidth}
\includegraphics[width=0.3\textwidth]{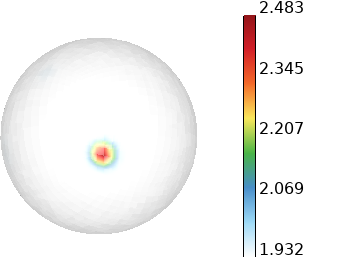}
\includegraphics[width=0.3\textwidth]{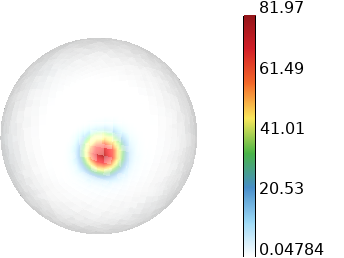}\hspace{0.2\textwidth}
\includegraphics[width=0.3\textwidth]{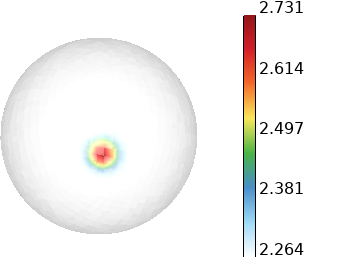}
\caption[Additional plots of $(u,v)$ with a random initial condition and $K=200000.0$]{Additional plots of $u$ (left) and $v$ (right) for a random initial condition and $K=200000.0$. The snapshots correspond to the time instants $t=1, 3$ (from top to bottom).}
\label{abb:Cgu2}
\end{figure}

\section{Conclusion}
In this paper we have studied a bulk-surface reaction-diffusion system of Gierer-Meinhardt type analytically and numerically. We were able to show the global in time existence of solutions to the fully coupled Problem \ref{sys:Ana}. For the proof we have used an operator splitting in a surface and a bulk component, and have applied the Schauder fixed point theorem. We have extended and carefully adapted previous results for the Gierer-Meinhardt system by Masuda and Takahashi \cite{MasTak1987} to obtain existence and a priori estimates for the surface system with fixed bulk contribution. For the bulk system with fixed input from the surface variables we used the theory of linear parabolic Robin boundary value problems \cite{LaSoUr1968}. 

In Section \ref{sec:3} we considered the reduced system for a well mixed bulk solution and discretized it using a tailor-made positivity-preserving finite element scheme. By simulating reaction-diffusion processes described by this system with different parameters and initial conditions we found localized steady-state multispike membrane-bound patterns and observed a competition between spikes. Depending on the parameter choices and the initial conditions steady states with one or two symmetrically distributed spikes were observed. This behaviour has also been found for the corresponding system on a one-dimensional boundary studied by Gomez et al in \cite{GoWaWe2018}. In addition they also observed oscillatory instabilities, which we have not seen in our simulations.

\bibliographystyle{spmpsci}
\bibliography{mybibfile}

\end{document}